\documentclass[hidelinks,onefignum,onetabnum]{siamart251216}

\usepackage{lipsum}
\usepackage{amsfonts}
\usepackage{graphicx}
\usepackage{epstopdf}
\usepackage{algorithmic}
\usepackage{amssymb}

\newcommand{\x}{\mathbf{x}}
\newcommand{\y}{\mathbf{y}}
\newcommand{\z}{\mathbf{z}}
\newcommand{\bd}{\mathbf{d}}
\newcommand{\bv}{\mathbf{v}}
\newcommand{\w}{\mathbf{w}}
\newcommand{\G}{\mathbf{G}}
\newcommand{\e}{\mathbf{e}}
\newcommand{\bu}{\mathbf{u}}

\DeclareMathOperator{\dist}{dist}
\DeclareMathOperator{\proj}{proj}

\ifpdf
  \DeclareGraphicsExtensions{.eps,.pdf,.png,.jpg}
\else
  \DeclareGraphicsExtensions{.eps}
\fi

\newsiamremark{remark}{Remark}
\newsiamremark{hypothesis}{Hypothesis}
\crefname{hypothesis}{Hypothesis}{Hypotheses}
\newsiamthm{claim}{Claim}
\newsiamremark{fact}{Fact}
\crefname{fact}{Fact}{Facts}

\headers{Stationarity Transfer and Repeatability}{S. Mei}

\title{Tight Stationarity Transfer and Repeatability
for Inexact Proximal-Gradient Methods}

\author{Shenghan Mei\thanks{Department of Mathematics,
University of North Carolina at Chapel Hill,
Chapel Hill, NC 27599, USA
(\email{shmei@unc.edu}).}}

\usepackage{amsopn}

\hypersetup{
  pdftitle={Tight Stationarity Transfer and Repeatability for Inexact Proximal-Gradient Methods},
  pdfauthor={Shenghan Mei}
}

\begin{document}

\maketitle


\begin{abstract}
Proximal-gradient methods are widely used for composite
optimization, but inexact proximal updates introduce a
residual in the proximal optimality condition, complicating
the analysis of how stationarity changes across a step.
We study this question for composite objectives consisting
of a smooth convex term and a possibly nonsmooth weakly
convex term, focusing on the relationship between the
directional-stationarity measures at the input and output
points of an inexact proximal-gradient step.
We first establish separate input and output bounds,
but their combination can yield a conservative estimate
because their worst-case residual directions may differ.
By exploiting the magnitude and orientation of this
shared residual, we establish tight joint bounds and
derive an exact one-step stationarity-transfer factor
that strictly improves the estimate from separate bounds
in the short-step regime.
We also identify the precise residual-tolerance threshold
for finite uniform stationarity transfer.
A further distinction emerges over multiple steps:
endpoint worst cases can repeat over arbitrarily many
iterations, whereas interior worst cases cannot persist
across two consecutive steps, yielding strict two-step
submultiplicativity.
These results also provide principled tolerance-selection
rules, descent estimates, and stationarity-complexity
guarantees for inexact proximal-gradient methods.
\end{abstract}



\section{Introduction}\label{sec:introduction}

Many optimization problems arising in machine learning,
signal processing, and inverse problems combine a smooth
data-fitting term with a possibly nonsmooth regularizer \cite{bach2011sparsity,combettes2005signal,beck2009fista}.
Proximal-gradient methods are well suited to such problems,
as they combine a gradient step for the smooth term with
a proximal update for the regularizer
\cite{beck2017firstorder,parikh2014proximal}.
For exact proximal-gradient steps, the proximal
optimality condition yields an upper bound on
the stationarity measure at the output point
\cite{chen2022proximal}.
In practice, however, the proximal subproblem may require
an iterative solver and thus be solved only approximately \cite{schmidt2011inexact,barre2023principled},
introducing a residual in its optimality condition.
Although this residual can be incorporated into an output
stationarity bound, such a bound alone does not describe
how stationarity at the output compares with that at
the input. This raises a natural question: how tightly
can the stationarity measures at the input and output
points of an inexact proximal-gradient step be related?

We address this question for composite objectives
consisting of a smooth convex term and a possibly
nonsmooth weakly convex term.
This setting includes nonconvex regularization models
while retaining first-order properties that facilitate
stationarity analysis \cite{davis2019stochastic}.
We use a directional-stationarity measure that vanishes
precisely at directionally stationary (d-stationary)
points \cite{ahn2017difference} and admits an equivalent
subdifferential representation in this setting
\cite{rockafellar1998variational,davis2019stochastic}.

Previous research on proximal-gradient methods has
addressed three closely related topics: stationarity,
inexactness, and worst-case performance.
For convex composite problems, relationships between
proximal-gradient mappings and subgradient-based
stationarity measures have been established
\cite{chen2022proximal}.
When proximal subproblems are solved approximately,
various inexactness criteria have been developed
\cite{barre2023principled,bellocruz2025relative,
lauga2026characterizations}.
Convergence guarantees have also been established for
inexact proximal-gradient and forward-backward splitting
methods
\cite{bonettini2020forwardbackward,
hamadouche2024sharper,gu2018inexact},
with further analyses for weakly
convex problems
\cite{khanh2025inexact,bednarczuk2026forwardbackward}.
A complementary line of research focuses on tight
worst-case guarantees, where performance-estimation
techniques have been used to analyze composite
first-order methods
\cite{taylor2017exact,taylor2018exactpg}.
Recent studies have also developed tight analyses of
the difference-of-convex algorithm (DCA) and its
connection to proximal-gradient descent under weak
curvature \cite{rotaru2025tight}, established
interpolation conditions for weakly convex functions
\cite{rubbens2026constraint}, and investigated
gradient descent with relative errors
\cite{vernimmen2026worstcase}.
These developments motivate our study of tight
input-output stationarity bounds for inexact
proximal-gradient steps and the repeatability of
their worst-case behavior over successive iterations.

To investigate this input-output relationship, we first
derive separate bounds for the two stationarity measures.
However, even when these bounds are individually tight,
their combination may be conservative because the
corresponding worst cases can require different
orientations of the same proximal residual.
To overcome this limitation, we retain the magnitude
and orientation of the shared residual and establish
a tight upper bound relating the normalized input
and output stationarity measures.
Using explicit two-dimensional quadratic constructions,
we demonstrate the tightness of this bound and derive
an exact closed-form one-step stationarity-transfer
factor.
In the short-step regime, this factor strictly improves
the estimate obtained from separate bounds for positive
relative residual tolerances below a threshold that
we identify exactly. Beyond this threshold, no finite
uniform stationarity-transfer factor exists.

While the one-step analysis provides an exact
stationarity-transfer factor, an important question
remains: can the corresponding worst-case behavior
persist over successive iterations?
A tight one-step bound does not necessarily yield
a tight multi-step bound, since consecutive steps
must be compatible with the same objective.
We show that the repeatability of the one-step worst
case depends on the orientation of the proximal residual. 
When the one-step worst case occurs at an endpoint,
the same worst-case factor can be attained at every
iteration over arbitrarily many steps.
In contrast, interior worst-case behavior
cannot persist across two consecutive steps.
We establish this nonrepeatability using cyclic
monotonicity, showing that the one-step worst case is transient rather than persistent.

Our analysis also determines how much inexactness
can be tolerated in proximal updates.
In particular, we identify the largest relative
residual tolerance that guarantees a prescribed
stationarity-transfer factor.
Under suitable conditions, we further establish
objective descent along the inexact proximal-gradient
iterates and quantify the number of iterations
needed to obtain an approximately d-stationary point.

The remainder of the paper is organized as follows.
Section~\ref{sec:setting} introduces the problem setting
and the directional-stationarity measure, and describes
inexact proximal-gradient steps through their
proximal optimality residuals.
Section~\ref{sec:stationarity-bounds} develops separate and
joint input-output stationarity bounds.
Section~\ref{sec:sharp-one-step} establishes a tight
input-output bound and derives the exact one-step
stationarity-transfer factor.
Section~\ref{sec:repeatability} investigates the
repeatability of worst-case behavior over multiple steps.
Finally, Section~\ref{sec:accuracy-descent} presents
criteria for selecting relative residual tolerances,
descent estimates, and stationarity-convergence
and complexity guarantees.

\section{Problem setting and preliminaries}\label{sec:setting}

In this section, we introduce the weakly convex composite setting,
the directional-stationarity measure used throughout the paper, and the residual model for inexact proximal-gradient steps.

We consider the problem
\begin{align}\label{prob:fplusg}
    \min_{\x \in \mathbb{R}^n}
    \ \zeta(\x) \triangleq f(\x) + g(\x),
\end{align}
where $f:\mathbb{R}^n\to\mathbb{R}$ is convex and continuously
differentiable with an $L$-Lipschitz continuous gradient, i.e.,
$
    \|\nabla f(\x)-\nabla f(\y)\|_2
    \le L\|\x-\y\|_2,
$
$\forall\,\x,\y\in\mathbb{R}^n$,
and $g:\mathbb{R}^n\to\mathbb{R}$ is weakly convex and possibly nonsmooth. The notion of weak convexity is defined in Definition \ref{def:weaklycvx}.

\begin{definition}[Weak convexity]\label{def:weaklycvx}
    A function $g:\mathbb{R}^n\to\mathbb{R}$ is called
    $\rho$-weakly convex for some $\rho\ge 0$ if
    $$
        \widehat{g}(\x)
        \triangleq
        g(\x)+\frac{\rho}{2}\|\x\|^2
    $$
    is convex. The constant $\rho$ is referred to as a weak convexity
    modulus of $g$.
\end{definition}

We adopt the above conditions on $f$ and $g$
as standing assumptions throughout the paper.

\subsection{Directional-stationarity measure}

For $\x,\bd\in\mathbb{R}^n$, the directional derivative of $\zeta$
at $\x$ along the direction $\bd$ is defined by
$$
    \zeta'(\x;\bd)
    \triangleq
    \lim_{\tau\downarrow0}
    \frac{\zeta(\x+\tau\bd)-\zeta(\x)}{\tau},
$$
whenever the limit exists. $\zeta$ is
directionally differentiable, since $g$ is weakly convex \cite[Proposition~2]{mei2026admm}. This motivates the following stationarity notion.

\begin{definition}[Directional stationarity]\label{def:Dstat}
A point $\x^*\in\mathbb{R}^n$ is called directionally stationary,
or d-stationary, for \eqref{prob:fplusg} if
$
    \zeta'(\x^*;\bd)\geq 0,
    \
    \forall\,\bd\in\mathbb{R}^n.
$
\end{definition}

We next quantify the extent to which a point fails to be
d-stationary.
\begin{definition}[Directional-stationarity measure]\label{def:d-gap}
For $\x\in\mathbb{R}^n$, the directional-stationarity measure is defined by
\begin{align}\label{eq:d-gap}
    \mathcal{D}_{\zeta}(\x)
    \triangleq
    -\min_{\|\bd\|_2\leq 1}\zeta'(\x;\bd).
\end{align}
\end{definition}
Since $\zeta'(\x;\mathbf{0})=0$, we have
$\mathcal{D}_{\zeta}(\x)\geq0$. Moreover,
$\mathcal{D}_{\zeta}(\x)=0$ if and only if $\x$ is d-stationary.

To connect the directional-stationarity measure with the convex
subdifferential of $\widehat{g}$, let $\partial\widehat{g}$ denote the convex
subdifferential of $\widehat{g}$ \cite{bauschke2017convex}, and define the stationarity set
of $\zeta$ at $\x$ by
$$
    \mathcal{S}_{\zeta}(\x)
    \triangleq
    \nabla f(\x)+\partial\widehat{g}(\x)-\rho\x.
$$
The following result shows that the directional-stationarity measure
is exactly the distance from the origin to this set.

\begin{lemma}[Stationarity-set representation]\label{lem:d-gap}
For every $\x\in\mathbb{R}^n$,
\begin{align}\label{eq:d-gap-rep}
    \mathcal{D}_{\zeta}(\x)
    =\dist \bigl(\mathbf{0},\mathcal{S}_{\zeta}(\x)\bigr).
\end{align}
Consequently, $\x$ is d-stationary if and only if
$\mathbf{0}\in\mathcal{S}_{\zeta}(\x)$.
\end{lemma}

\begin{proof}
Since $\widehat{g}$ is convex, its directional
derivative is the support function of its convex subdifferential
\cite[Theorem~8.30]{rockafellar1998variational}, i.e.,
$
    \widehat{g}'(\x;\bd)
    =
    \sup_{\bv\in\partial\widehat{g}(\x)}
    \langle \bv,\bd\rangle.
$
Using
$
    g(\x)
    =
    \widehat{g}(\x)-\frac{\rho}{2}\|\x\|_2^2,
$
we obtain
$
    g'(\x;\bd)
    =
    \widehat{g}'(\x;\bd)-\rho\langle\x,\bd\rangle.
$
Since $f$ is continuously differentiable,
$
    f'(\x;\bd)=\langle\nabla f(\x),\bd\rangle.
$
Therefore,
\begin{align*}
    \zeta'(\x;\bd)
    &=
    \langle\nabla f(\x),\bd\rangle
    +\widehat{g}'(\x;\bd)
    -\rho\langle\x,\bd\rangle \\
    &=
    \sup_{\bv\in\partial\widehat{g}(\x)}
    \left\langle
        \nabla f(\x)+\bv-\rho\x,\bd
    \right\rangle \\
    &=
    \sup_{\w\in\mathcal{S}_{\zeta}(\x)}
    \langle\w,\bd\rangle.
\end{align*}

Let
$
    \w^\star
    =
    \proj_{\mathcal{S}_{\zeta}(\x)}(\mathbf{0}).
$
Since $\mathcal{S}_{\zeta}(\x)$ is nonempty, closed, and convex, 
$
\|\w^\star\|_2=\dist\bigl(\mathbf{0},\mathcal{S}_{\zeta}(\x)\bigr).
$
For any $\|\bd\|_2\leq1$,
$$
    \zeta'(\x;\bd)
    =
    \sup_{\w\in\mathcal{S}_{\zeta}(\x)}
    \langle\w,\bd\rangle
    \geq
    \langle\w^\star,\bd\rangle
    \geq
    -\|\w^\star\|_2.
$$
Hence
$
    \min_{\|\bd\|_2\leq1}\zeta'(\x;\bd)
    \geq
    -\|\w^\star\|_2.
$

If $\w^\star\neq\mathbf{0}$, define
$
    \bd^\star
    =
    -\frac{\w^\star}{\|\w^\star\|_2}.
$
By the projection inequality,
$
    \langle\w-\w^\star,\w^\star\rangle\geq0,    
$
for every $\w\in\mathcal{S}_{\zeta}(\x)$,
and therefore
$
    \langle\w,\bd^\star\rangle
    \leq
    -\|\w^\star\|_2,
    \
    \forall\,\w\in\mathcal{S}_{\zeta}(\x).
$
Since $\w^\star\in\mathcal{S}_{\zeta}(\x)$, equality is attained at
$\w=\w^\star$, and thus
$
    \zeta'(\x;\bd^\star)
    =
    -\|\w^\star\|_2.
$
Consequently,
$
    \min_{\|\bd\|_2\leq1}\zeta'(\x;\bd)
    =
    -\|\w^\star\|_2.
$
If $\w^\star=\mathbf{0}$, then $\mathbf{0}\in\mathcal{S}_{\zeta}(\x)$,
so $\zeta'(\x;\bd)\geq0$ for all $\bd$, while
$\zeta'(\x;\mathbf{0})=0$; hence both sides of
\eqref{eq:d-gap-rep} are zero. Therefore
$
    \mathcal{D}_{\zeta}(\x)
    =
    \|\w^\star\|_2
    =
    \operatorname{dist}
    \bigl(\mathbf{0},\mathcal{S}_{\zeta}(\x)\bigr).
$

Finally,
$\mathcal{D}_{\zeta}(\x)=0$ if and only if
$\mathbf{0}\in\mathcal{S}_{\zeta}(\x)$, which is equivalent to
$\x$ being d-stationary.
\end{proof}

\subsection{Inexact proximal-gradient steps}

We next introduce the residual model used to quantify inexactness in a
proximal-gradient step. Given a current point $\x\in\mathbb{R}^n$ and
a stepsize $\alpha>0$, the exact proximal-gradient update is defined by
\cite{bauschke2017convex}
$$
    \z
    \in
    \arg\min_{\y\in\mathbb{R}^n}
    \left\{
        g(\y)
        +
        \langle\nabla f(\x),\y-\x\rangle
        +
        \frac{1}{2\alpha}\|\y-\x\|_2^2
    \right\}.
$$
When $\alpha\rho<1$, the above subproblem is strongly convex and hence
admits a unique minimizer.

For the input point $\x$ and a candidate output point $\z$, define the proximal-gradient displacement
and the proximal optimality-residual set by
\begin{align}\label{eq:PG-residual}
    \G
    &=
    \frac{\x-\z}{\alpha},
    \qquad
    \mathcal{R}_{\alpha}(\x,\z)
    \triangleq
    \partial\widehat{g}(\z)
    +
    \nabla f(\x)
    -
    \rho\z
    -
    \G.
\end{align}
A vector $\e\in\mathcal{R}_{\alpha}(\x,\z)$ is called a residual
vector. Equivalently, there exists
$\bv_{\z}\in\partial\widehat{g}(\z)$ such that
$
    \e
    =
    \bv_{\z}
    +
    \nabla f(\x)
    -
    \rho\z
    -
    \G.
$
Indeed, using
$
g(\y)=\widehat{g}(\y)-\frac{\rho}{2}\|\y\|_2^2,
$
the first-order optimality condition of the proximal subproblem is
precisely
$
\mathbf{0}\in\mathcal{R}_{\alpha}(\x,\z).
$
Thus, when $\alpha\rho<1$, this condition uniquely characterizes the
exact proximal-gradient update, while a nonzero residual vector
quantifies the violation of the exact proximal optimality condition.

A central inexactness criterion in our analysis is the relative
residual condition
\begin{align}\label{eq:relative-residual}
    \|\e\|_2
    \leq
    \eta\|\G\|_2,
\end{align}
where $\eta\geq0$ is the prescribed relative tolerance. Whenever
$\G\neq\mathbf{0}$, define
\begin{align}\label{eq:normalized-residual}
    \bu
    =
    \frac{\G}{\|\G\|_2},
    \qquad
    q
    =
    \frac{\|\e\|_2}{\|\G\|_2},
    \qquad
    s
    = \frac{\langle\e,\G\rangle}{\|\G\|_2^2}.
\end{align}
Here, $q$ measures the relative magnitude of the residual, while $s$
is its signed component along the computed displacement. By the
Cauchy-Schwarz inequality,
$
    |s|\leq q,
$
and \eqref{eq:relative-residual} implies $q\leq\eta$.

For convenience, we use the following quantities
throughout the subsequent analysis:
$$
a_{\alpha}=1-\alpha\rho,
\qquad
\theta_{\alpha}=\frac{\alpha L}{2},
\qquad
c_{\alpha}=1-\frac{\alpha L}{2}.
$$
We also use
$
[r]_+ = \max\{r,0\}
$
for the positive part of a scalar $r$.

\section{Input-output stationarity bounds}
\label{sec:stationarity-bounds}

In this section, we establish stationarity bounds for
inexact proximal-gradient steps. 
We first establish separate input and output bounds
and then derive tighter joint bounds by retaining
the magnitude and orientation of the shared residual.

\subsection{Separate input and output bounds}
\label{sec:separate-bounds}

For exact proximal-gradient steps in convex composite
optimization, bounds relating proximal-gradient mappings
to stationarity measures are known
\cite{chen2022proximal}. We extend these separate input and output bounds to weakly convex objectives and inexact proximal solves.
For convenience, set
$
    \kappa_{\alpha}=\max\{1,\alpha L-1\}.
$

\begin{proposition}[Input and output stationarity bounds]
\label{prop:marginal}
Let $\e\in\mathcal{R}_{\alpha}(\x,\z)$. Then
\begin{equation}\label{eq:separate-bounds}
\begin{aligned}
    \mathcal{D}_{\zeta}(\z)
    &\leq
    \kappa_{\alpha}\|\G\|_2+\|\e\|_2,\\
    a_{\alpha}\|\G\|_2
    &\leq
    \mathcal{D}_{\zeta}(\x)+\|\e\|_2.
\end{aligned}
\end{equation}
Consequently, if
$
\|\e\|_2\leq\eta\|\G\|_2,
$
then
\begin{equation}\label{eq:separate-relative-bounds}
\begin{aligned}
    \mathcal{D}_{\zeta}(\z)
    &\leq
    (\kappa_{\alpha}+\eta)\|\G\|_2,\\
    \mathcal{D}_{\zeta}(\x)
    &\geq
    [a_{\alpha}-\eta]_+\|\G\|_2.
\end{aligned}
\end{equation}
Both coefficients in \eqref{eq:separate-relative-bounds} are
individually best possible over the stated problem class.
\end{proposition}

\begin{proof}
Let
$
\Delta=\nabla f(\z)-\nabla f(\x).
$
By the definition of a residual vector, there exists
$\bv_{\z}\in\partial\widehat{g}(\z)$ such that
$
    \e
    =
    \bv_{\z}
    +
    \nabla f(\x)
    -
    \rho\z
    -
    \G.
$
Hence
$
    \G+\Delta+\e
    =
    \nabla f(\z)+\bv_{\z}-\rho\z
    \in
    \mathcal{S}_{\zeta}(\z).
$
Lemma~\ref{lem:d-gap} therefore gives
$$
    \mathcal{D}_{\zeta}(\z)
    \leq
    \|\G+\Delta+\e\|_2
    \leq
    \|\G+\Delta\|_2+\|\e\|_2.
$$

We next bound $\|\G+\Delta\|_2$. If $L=0$, then $\Delta=\mathbf{0}$
and the claim is immediate. Suppose $L>0$. Since $f$ is convex with an
$L$-Lipschitz continuous gradient, $\nabla f$ is $\frac{1}{L}$-cocoercive \cite[Corollary 18.17]{bauschke2017convex}.
Using $\z-\x=-\alpha\G$, we obtain
$
    \langle \G,\Delta\rangle
    \leq
    -\frac{1}{\alpha L}\|\Delta\|_2^2.
$
Therefore,
$$
    \|\G+\Delta\|_2^2
    \leq
    \|\G\|_2^2
    +
    \left(1-\frac{2}{\alpha L}\right)\|\Delta\|_2^2.
$$
If $\alpha L<2$, this gives
$
\|\G+\Delta\|_2\leq\|\G\|_2.
$
If $\alpha L\geq2$, the Lipschitz continuity of $\nabla f$ yields
$
\|\Delta\|_2\leq\alpha L\|\G\|_2,
$
and hence
$
    \|\G+\Delta\|_2
    \leq
    (\alpha L-1)\|\G\|_2.
$
Thus, in both cases,
$
    \|\G+\Delta\|_2
    \leq
    \kappa_{\alpha}\|\G\|_2,
$
which proves the output estimate in
\eqref{eq:separate-bounds}.

For the input estimate, let
$
    \mathbf{w}_{\x}
    =
    \proj_{\mathcal{S}_{\zeta}(\x)}(\mathbf{0}),
$
so that
$
\|\mathbf{w}_{\x}\|_2=\mathcal{D}_{\zeta}(\x).
$
There exists
$\bv_{\x}\in\partial\widehat{g}(\x)$ such that
$
    \mathbf{w}_{\x}
    =
    \nabla f(\x)+\bv_{\x}-\rho\x.
$
Since $\widehat{g}$ is convex, the monotonicity of its subdifferential
gives
$
    \langle
        \bv_{\z}-\bv_{\x},
        \z-\x
    \rangle
    \geq0.
$
From the definitions of the residual vector $\e$ and
$\mathbf{w}_{\x}$, we have
$
\bv_{\z}-\bv_{\x}
=
\e+\G-\mathbf{w}_{\x}
+\rho(\z-\x).
$
Substituting this identity together with
$
\z-\x=-\alpha\G
$
into the monotonicity inequality above yields
$$
\langle
\mathbf{w}_{\x}-\e,
\G
\rangle
\geq
(1-\alpha\rho)\|\G\|_2^2
=
a_{\alpha}\|\G\|_2^2.
$$
If $\G=\mathbf{0}$, the input estimate is immediate. Otherwise,
Cauchy-Schwarz gives
$$
    a_{\alpha}\|\G\|_2
    \leq
    \|\mathbf{w}_{\x}\|_2+\|\e\|_2
    =
    \mathcal{D}_{\zeta}(\x)+\|\e\|_2,
$$
which proves the input estimate in \eqref{eq:separate-bounds}.
Applying
$
\|\e\|_2\leq\eta\|\G\|_2
$
to both inequalities in \eqref{eq:separate-bounds} yields
\eqref{eq:separate-relative-bounds}.

It remains to show that the two coefficients in
\eqref{eq:separate-relative-bounds} are individually best possible. For this purpose, fix any $\G\neq\mathbf{0}$. For the output bound, if
$\alpha L<2$, take $g=0$ and an affine $f$ with
$
\nabla f=(1+\eta)\G;
$
then $\e=\eta\G$ and
$
\mathcal{D}_{\zeta}(\z)=(1+\eta)\|\G\|_2.
$
If $\alpha L\geq2$, take
$
f(\y)=\frac{L}{2}\|\y\|_2^2
$
and $g=0$. Choose any nonzero $\G$ and set
$
\x=\frac{1-\eta}{L}\G
$
and
$
\z=\x-\alpha\G.
$
Then $\e=-\eta\G$ and
$
\mathcal{D}_{\zeta}(\z)
=
(\alpha L-1+\eta)\|\G\|_2.
$

For the input bound, suppose first that $a_{\alpha}>\eta$.
If $\rho>0$, choose any $\G\neq\mathbf{0}$ and set
$
f=0
$
and
$
g(\y)=-\frac{\rho}{2}\|\y\|_2^2.
$
Let
$
\x=-\frac{a_{\alpha}-\eta}{\rho}\G
$
and
$
\z=\x-\alpha\G.
$
Since $\widehat{g}=0$, we have
$
\mathcal{S}_{\zeta}(\x)=\{-\rho\x\}
$
and
$
\e=-\rho\z-\G=-\eta\G.
$
Therefore,
$
\mathcal{D}_{\zeta}(\x)
=
(a_{\alpha}-\eta)\|\G\|_2,
$
so the positive input coefficient in
\eqref{eq:separate-relative-bounds} is attained.
If $\rho=0$, then $a_{\alpha}=1$. Taking $g=0$ and the affine function
$
f(\y)=\langle(1-\eta)\G,\y\rangle
$
gives
$
\e=-\eta\G
$
and
$
\mathcal{D}_{\zeta}(\x)
=
(1-\eta)\|\G\|_2
=
(a_{\alpha}-\eta)\|\G\|_2.
$

Finally, suppose that $a_{\alpha}\leq\eta$. Choose any
$\G\neq\mathbf{0}$ and set
$
f=0,
$
$
\x=\mathbf{0},
$
and
$
\z=-\alpha\G.
$
Let
$
m=[-a_{\alpha}]_+/\alpha
$
and define
$
g(\y)=\frac{m-\rho}{2}\|\y\|_2^2.
$
Then
$
\widehat{g}(\y)=\frac{m}{2}\|\y\|_2^2
$
is convex, so $g$ is $\rho$-weakly convex. Moreover,
$
\mathcal{D}_{\zeta}(\x)=0
$
and
$
\e=-(a_{\alpha}+\alpha m)\G
=-[a_{\alpha}]_+\G.
$
Since $a_{\alpha}\leq\eta$, we have
$
[a_{\alpha}]_+\leq\eta,
$
and hence
$
\|\e\|_2\leq\eta\|\G\|_2.
$
Thus the zero coefficient in
$
[a_{\alpha}-\eta]_+
$
is also best possible.
\end{proof}

When $\eta<a_{\alpha}$, the input coefficient in
\eqref{eq:separate-relative-bounds} is positive. Combining the two
bounds in Proposition~\ref{prop:marginal} therefore gives
\begin{align}\label{eq:marginal-ratio}
    \mathcal{D}_{\zeta}(\z)
    \leq
    \frac{\kappa_{\alpha}+\eta}
         {a_{\alpha}-\eta}
    \mathcal{D}_{\zeta}(\x).
\end{align}
The two coefficients in \eqref{eq:separate-relative-bounds} are
individually best possible, but this does not imply that the resulting
stationarity-transfer bound is best possible. Indeed, the input and
output bounds are optimized separately, whereas in a single
proximal-gradient step they depend on the same residual vector $\e$.
The next subsection accounts for this coupling and derives a tighter
joint stationarity-transfer bound.

\subsection{Joint input and output bounds}
\label{sec:joint-geometry}

The separate bounds in Section \ref{sec:separate-bounds} use the residual only through its
norm and therefore do not retain its orientation relative to the
proximal-gradient displacement. We now keep this information through
the quantities $q$ and $s$ introduced in Section \ref{sec:setting}. In particular,
$q$ measures the relative residual magnitude, while $s$ measures its
signed component along $\G$. The key point is that both the input and output bounds depend on the
same residual vector $\e$. By retaining the residual quantities $q$
and $s$, the resulting joint bounds can be tighter than the separate
estimates. To derive the output bound, we first characterize
$\nabla f(\z)-\nabla f(\x)$ using cocoercivity.

\begin{lemma}[Cocoercivity ball]\label{lem:coco-ball}
Suppose $\G\neq\mathbf{0}$ and define
$
\mathbf{h}
=
-\frac{\nabla f(\z)-\nabla f(\x)}{\|\G\|_2}.
$
Then
$$
\|\mathbf{h}-\theta_{\alpha}\bu\|_2
\leq
\theta_{\alpha}.
$$
In particular, if $L=0$, then
$\mathbf{h}=\mathbf{0}$ and $\theta_{\alpha}=0$.
\end{lemma}

\begin{proof}
If $L=0$, then $\nabla f$ is constant, so
$\mathbf{h}=\mathbf{0}$ and the conclusion is immediate.
Suppose now that $L>0$. By the cocoercivity of $\nabla f$,
$$
\left\langle
\nabla f(\z)-\nabla f(\x),\z-\x
\right\rangle
\geq
\frac{1}{L}
\|\nabla f(\z)-\nabla f(\x)\|_2^2.
$$
Using $\z-\x=-\alpha\G$ and the definitions of
$\mathbf{h}$ and $\bu$ gives
$
\langle\mathbf{h},\bu\rangle
\geq
\frac{\|\mathbf{h}\|_2^2}{\alpha L}.
$
Since $\theta_{\alpha}=\frac{\alpha L}{2}$, we obtain
$$
\begin{aligned}
\|\mathbf{h}-\theta_{\alpha}\bu\|_2^2
&=
\|\mathbf{h}\|_2^2
-2\theta_{\alpha}\langle\mathbf{h},\bu\rangle
+\theta_{\alpha}^2 \\
&\leq
\|\mathbf{h}\|_2^2
-\frac{2\theta_{\alpha}}{\alpha L}\|\mathbf{h}\|_2^2
+\theta_{\alpha}^2
=
\theta_{\alpha}^2,
\end{aligned}
$$
which proves the result.
\end{proof}

The next theorem gives joint input and output stationarity bounds that
depend explicitly on the residual quantities $q$ and $s$.

\begin{theorem}[Joint input-output stationarity bounds]
\label{thm:joint-bounds}
Let $\e\in\mathcal{R}_{\alpha}(\x,\z)$ and suppose that
$\G\neq\mathbf{0}$. Then
\begin{align}
\mathcal{D}_{\zeta}(\x)
&\geq
[a_{\alpha}+s]_+\|\G\|_2,
\label{eq:joint-input}\\
\mathcal{D}_{\zeta}(\z)
&\leq
\left(
\theta_{\alpha}
+
\sqrt{
c_{\alpha}^2+2c_{\alpha}s+q^2
}
\right)\|\G\|_2.
\label{eq:joint-output}
\end{align}
\end{theorem}

\begin{proof}
For the input bound, the argument in the proof of
Proposition~\ref{prop:marginal} gives
$
\langle\w_{\x}-\e,\G\rangle
\geq
a_{\alpha}\|\G\|_2^2,
$
where
$
\w_{\x}
=
\operatorname{proj}_{\mathcal{S}_{\zeta}(\x)}(\mathbf{0})
$
and hence
$
\|\w_{\x}\|_2=\mathcal{D}_{\zeta}(\x).
$
Since
$
\langle\e,\G\rangle=s\|\G\|_2^2,
$
we obtain
$$
\langle\w_{\x},\G\rangle
\geq
(a_{\alpha}+s)\|\G\|_2^2.
$$
If $a_{\alpha}+s>0$, the Cauchy-Schwarz inequality gives
$
\mathcal{D}_{\zeta}(\x)
\geq
(a_{\alpha}+s)\|\G\|_2.
$
If $a_{\alpha}+s\leq0$, the result follows from the nonnegativity of
$\mathcal{D}_{\zeta}(\x)$. This proves \eqref{eq:joint-input}.

For the output bound, let $\mathbf{h}$ be as in
Lemma \ref{lem:coco-ball}.
As in the proof of Proposition \ref{prop:marginal},
$
\mathcal{D}_{\zeta}(\z)
\leq
\left\|
\G+\nabla f(\z)-\nabla f(\x)+\e
\right\|_2.
$
Hence
$
\frac{\mathcal{D}_{\zeta}(\z)}{\|\G\|_2}
\leq
\left\|
\bu-\mathbf{h}
+\frac{\e}{\|\G\|_2}
\right\|_2.
$
By Lemma~\ref{lem:coco-ball},
$
\|\mathbf{h}-\theta_{\alpha}\bu\|_2\leq\theta_{\alpha}.
$
Using $c_{\alpha}=1-\theta_{\alpha}$ and the triangle inequality,
$$
\begin{aligned}
\left\|
\bu-\mathbf{h}
+\frac{\e}{\|\G\|_2}
\right\|_2
&=
\left\|
c_{\alpha}\bu
+\frac{\e}{\|\G\|_2}
-
(\mathbf{h}-\theta_{\alpha}\bu)
\right\|_2\\
&\leq
\left\|
c_{\alpha}\bu
+\frac{\e}{\|\G\|_2}
\right\|_2
+\theta_{\alpha}.
\end{aligned}
$$
Moreover, since $\|\bu\|_2=1$, the definitions of $q$ and $s$ give
$$
\begin{aligned}
\left\|
c_{\alpha}\bu
+\frac{\e}{\|\G\|_2}
\right\|_2^2
&=
c_{\alpha}^2
+
2c_{\alpha}
\frac{\langle\e,\G\rangle}{\|\G\|_2^2}
+
\frac{\|\e\|_2^2}{\|\G\|_2^2}\\
&=
c_{\alpha}^2+2c_{\alpha}s+q^2.
\end{aligned}
$$
Combining these inequalities yields \eqref{eq:joint-output}.
\end{proof}

Under the relative residual condition
$\|\e\|_2\leq\eta\|\G\|_2$, we have $|s|\leq q\leq\eta$.
Consequently,
$
[a_{\alpha}+s]_+\geq[a_{\alpha}-\eta]_+
$
and
$
\theta_{\alpha}
+
\sqrt{c_{\alpha}^2+2c_{\alpha}s+q^2}
\leq
\kappa_{\alpha}+\eta.
$
Thus the bounds in Theorem~\ref{thm:joint-bounds} are at least as
tight as the corresponding bounds in Proposition~\ref{prop:marginal}. These inequalities can be strict
because the joint bounds retain the actual residual quantities $q$
and $s$.
Combining the two estimates in Theorem~\ref{thm:joint-bounds} yields
a stationarity-transfer bound for the current step.

\begin{corollary}[Residual-based stationarity transfer]
\label{cor:residual-transfer}
Under the assumptions of Theorem~\ref{thm:joint-bounds}, suppose that
$a_{\alpha}+s>0$. Define
$
\Theta_{\alpha}(q,s)
=
\frac{
\theta_{\alpha}
+
\sqrt{c_{\alpha}^2+2c_{\alpha}s+q^2}
}{
a_{\alpha}+s
}.
$
Then
$$
\mathcal{D}_{\zeta}(\z)
\leq
\Theta_{\alpha}(q,s)\,
\mathcal{D}_{\zeta}(\x).
$$
Consequently, for any $\gamma>0$, if
$
\Theta_{\alpha}(q,s)\leq\gamma,
$
then
$
\mathcal{D}_{\zeta}(\z)
\leq
\gamma\mathcal{D}_{\zeta}(\x).
$
\end{corollary}

\begin{proof}
Since $a_{\alpha}+s>0$, combining
\eqref{eq:joint-input} and \eqref{eq:joint-output} gives
$$
\mathcal{D}_{\zeta}(\z)
\leq
\frac{
\theta_{\alpha}
+
\sqrt{c_{\alpha}^2+2c_{\alpha}s+q^2}
}{
a_{\alpha}+s
}
\mathcal{D}_{\zeta}(\x)
=
\Theta_{\alpha}(q,s)\mathcal{D}_{\zeta}(\x).
$$
The final statement follows immediately if
$\Theta_{\alpha}(q,s)\leq\gamma$.
\end{proof}

The bounds in Theorem~\ref{thm:joint-bounds} still depend on the
computed residual quantities $q$ and $s$. Under the relative
residual condition, these quantities satisfy $|s|\leq q\leq\eta$.
The next section removes the dependence on the computed residual
quantities $q$ and $s$ and derives bounds depending only on the
relative residual tolerance $\eta$.

\section{Tight input-output bounds and exact one-step factor}
\label{sec:sharp-one-step}

We now impose the relative residual condition
$\|\e\|_2\leq\eta\|\G\|_2$. Since $|s|\leq q\leq\eta$, we optimize
the joint bounds in Theorem~\ref{thm:joint-bounds} over all admissible
pairs $(q,s)$. We first consider the regime
$0\leq\eta<a_{\alpha}$, for which $a_{\alpha}+s>0$ for every
admissible $s$.

\subsection{Tight input-output bound}
\label{subsec:sharp-frontier}

For $\G\neq\mathbf{0}$, define the normalized stationarity measures
$
X
=
\frac{\mathcal{D}_{\zeta}(\x)}{\|\G\|_2}
\text{ and }
Z
=
\frac{\mathcal{D}_{\zeta}(\z)}{\|\G\|_2}.
$
We establish a tight upper bound on $Z$ in terms of $X$.

\begin{theorem}[Tight input-output bound]
\label{thm:sharp-frontier}
Suppose that $0\leq\eta<a_{\alpha}$ and
$\|\e\|_2\leq\eta\|\G\|_2$. Then
$
X\geq a_{\alpha}-\eta.
$
Moreover, the following statements hold.

\emph{(i) Short-step regime.}
If $\alpha L<2$, then
\begin{align}
Z
\leq
\theta_{\alpha}
+
\sqrt{
c_{\alpha}^2
+
2c_{\alpha}\min\{\eta,X-a_{\alpha}\}
+
\eta^2
}.
\label{eq:frontier-short}
\end{align}

\emph{(ii) Long-step regime.}
If $\alpha L\geq2$, then
$
Z\leq\alpha L-1+\eta.
$

\end{theorem}

\begin{proof}
Let $q$ and $s$ be the residual quantities from Section \ref{sec:setting}. Since
$
|s|\leq q\leq\eta,
$
the input bound \eqref{eq:joint-input} gives
$
X\geq a_{\alpha}+s\geq a_{\alpha}-\eta.
$
The output bound \eqref{eq:joint-output} gives
$$
Z
\leq
\theta_{\alpha}
+
\sqrt{
c_{\alpha}^2+2c_{\alpha}s+\eta^2
}.
$$

Suppose first that $\alpha L<2$, so that $c_{\alpha}>0$. The right-hand
side is increasing in $s$. Moreover,
$
X\geq a_{\alpha}+s
$
implies
$
s\leq X-a_{\alpha},
$
while $s\leq\eta$. Hence
$
s\leq\min\{\eta,X-a_{\alpha}\},
$
which proves \eqref{eq:frontier-short}.

If $\alpha L\geq2$, then $c_{\alpha}\leq0$, so the output bound is
nonincreasing in $s$. Its largest value is therefore obtained at
$s=-\eta$, giving
$$
\theta_{\alpha}
+
\sqrt{c_{\alpha}^2-2c_{\alpha}\eta+\eta^2}
=
\alpha L-1+\eta,
$$
which proves the long-step bound and completes the proof.
\end{proof}

To show that the upper bound in Theorem~\ref{thm:sharp-frontier} is tight, we construct a family of two-dimensional quadratic examples
whose normalized stationarity pair $(X,Z)$ approach its boundary. The same construction will also be used later in the
threshold analysis.

\begin{proposition}[Tightness construction]
\label{prop:frontier-construction}
For any $\eta\geq0$, $s\in[-\eta,\eta]$, and
$\delta>0$ satisfying $a_{\alpha}+s+\delta>0$, there exist quadratic functions
$f,g_{\delta}:\mathbb{R}^2\to\mathbb{R}$ and an inexact
proximal-gradient step satisfying
$
\|\e\|_2=\eta\|\G\|_2
$
and
$
\frac{\langle\e,\G\rangle}{\|\G\|_2^2}=s,
$
for which
$$
X=a_{\alpha}+s+\delta,
\qquad
Z=
\theta_{\alpha}
+
\sqrt{c_{\alpha}^2+2c_{\alpha}s+\eta^2}.
$$
\end{proposition}

\begin{proof}
Let $\bu,\bv\in\mathbb{R}^2$ be orthonormal and set
$
r=\sqrt{\eta^2-s^2}
$
and
$
\boldsymbol{\xi}=s\bu+r\bv.
$
Set
$
\x=\mathbf{0},
$
$
\z=-\alpha\bu,
$
$
\G=\bu,
$
and
$
\e=\boldsymbol{\xi}.
$
Then
$
\|\G\|_2=1,
$
$
\|\e\|_2=\eta,
$
and
$
\langle\e,\G\rangle=s.
$

We first construct $f$. Let
$
Q=\|c_{\alpha}\bu+\boldsymbol{\xi}\|_2
$
and define
$$
\mathbf{h}
=
\begin{cases}
\displaystyle
\theta_{\alpha}\bu
-
\theta_{\alpha}
\frac{c_{\alpha}\bu+\boldsymbol{\xi}}{Q},
& Q>0,\\[3mm]
\mathbf{0},
& Q=0.
\end{cases}
$$
If $\mathbf{h}\neq\mathbf{0}$, let
$
A=L\frac{\mathbf{h}\mathbf{h}^{\top}}{\|\mathbf{h}\|_2^2};
$
otherwise, let $A=0$. Define
$
f(\w)=\frac12\langle\w,A\w\rangle.
$

When $Q>0$, the choice of $\mathbf{h}$ gives
$
\|\mathbf{h}-\theta_{\alpha}\bu\|_2=\theta_{\alpha},
$
and hence
$
\|\mathbf{h}\|_2^2
=
\alpha L\langle\mathbf{h},\bu\rangle.
$
Thus, when $\mathbf{h}\neq\mathbf{0}$,
$
A\bu
=
L\frac{\mathbf{h}\langle\mathbf{h},\bu\rangle}
        {\|\mathbf{h}\|_2^2}
=
\frac{\mathbf{h}}{\alpha}.
$
Moreover, $A$ is a positive semidefinite matrix with
$\|A\|_2=L$. If $\mathbf{h}=\mathbf{0}$, then $A=0$, so the same
properties hold with $\|A\|_2\leq L$. Therefore $f$ is convex with an
$L$-Lipschitz continuous gradient, and
$
\nabla f(\x)=\mathbf{0}
$
and
$
\nabla f(\z)=-\mathbf{h}.
$

We next construct $g_{\delta}$. For $\delta>0$, let
$
M_{\delta}
=
(\alpha\delta)^{-1}
(\delta\bu-r\bv)(\delta\bu-r\bv)^{\top},
$
and define
$$
\begin{aligned}
\widehat{g}_{\delta}(\w)
&=
\frac12\langle\w,M_{\delta}\w\rangle
+
(a_{\alpha}+s+\delta)\langle\bu,\w\rangle,\\
g_{\delta}(\w)
&=
\widehat{g}_{\delta}(\w)
-\frac{\rho}{2}\|\w\|_2^2.
\end{aligned}
$$
Since $M_{\delta}\succeq0$, the function
$\widehat{g}_{\delta}$ is convex. Hence $g_{\delta}$ is
$\rho$-weakly convex.

At the input point $\x=\mathbf{0}$,
$
\mathcal{S}_{\zeta}(\x)
=
\{(a_{\alpha}+s+\delta)\bu\}.
$
Since $a_{\alpha}+s+\delta>0$ and $\|\G\|_2=1$, we obtain
$
X=a_{\alpha}+s+\delta.
$

We next verify that the constructed functions generate the prescribed
residual $\e=\boldsymbol{\xi}$.
Since
$
M_{\delta}\bu=\frac{\delta\bu-r\bv}{\alpha},
$
using $\z=-\alpha\bu$ gives
$
\nabla\widehat{g}_{\delta}(\z)
=
(a_{\alpha}+s)\bu+r\bv.
$
Together with
$
\nabla f(\x)=\mathbf{0},
$
$
-\rho\z=\alpha\rho\bu,
$
and
$
a_{\alpha}=1-\alpha\rho,
$
this yields
$$
\nabla\widehat{g}_{\delta}(\z)
+\nabla f(\x)
-\rho\z-\G
=
s\bu+r\bv
=
\boldsymbol{\xi}
=
\e.
$$
Hence
$
\e\in\mathcal{R}_{\alpha}(\x,\z)
$
and
$
\|\e\|_2=\eta\|\G\|_2.
$

Finally, at the output point,
$
\mathcal{S}_{\zeta}(\z)
=
\{\bu+\boldsymbol{\xi}-\mathbf{h}\}.
$
By the choice of $\mathbf{h}$,
$
\|\bu+\boldsymbol{\xi}-\mathbf{h}\|_2
=
\theta_{\alpha}
+
\|c_{\alpha}\bu+\boldsymbol{\xi}\|_2.
$
Since
$
\|c_{\alpha}\bu+\boldsymbol{\xi}\|_2^2
=
c_{\alpha}^2+2c_{\alpha}s+\eta^2,
$
Lemma~\ref{lem:d-gap} yields
$$
Z
=
\theta_{\alpha}
+
\sqrt{c_{\alpha}^2+2c_{\alpha}s+\eta^2},
$$
which gives the result.
\end{proof}

The construction is essentially two-dimensional and extends to every
fixed dimension $n\geq2$ by setting the quadratic forms to zero on the
orthogonal complement of $\text{span}\{\bu,\bv\}$. When $0\leq\eta<a_{\alpha}$, we have
$a_{\alpha}+s>0$ for every $s\in[-\eta,\eta]$.
Thus, letting $\delta\downarrow0$ gives
$X\to a_{\alpha}+s$ while $Z$ remains unchanged,
establishing that the upper bound in Theorem~\ref{thm:sharp-frontier} is tight in every fixed
dimension $n\geq2$. At the endpoints $|s|=\eta$, the orthogonal
residual component vanishes, and the corresponding boundary values
are attained exactly by a one-dimensional construction.

\begin{figure}[htbp]
\centering
\includegraphics[width=0.65\linewidth]
{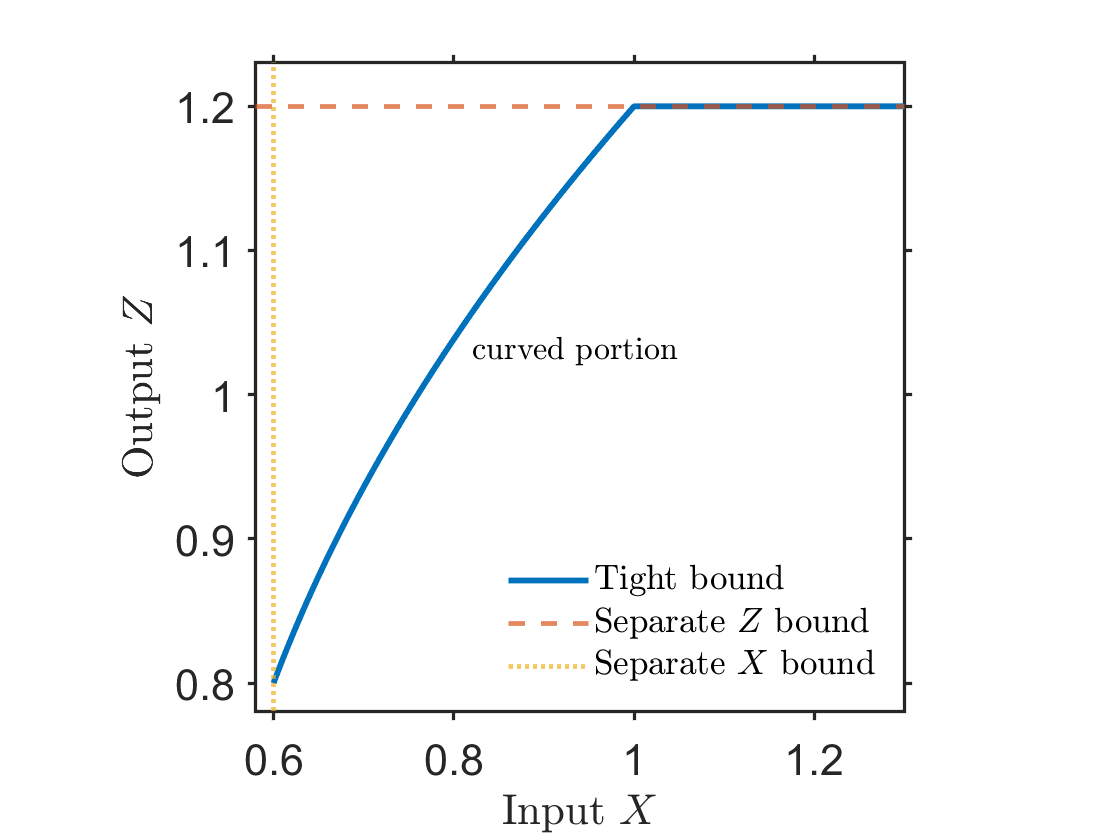}
\caption{
Tight upper bound on $Z$ in terms of $X$ for
$\alpha L=1$, $\alpha\rho=0.2$, and $\eta=0.2$.
The solid curve shows the tight upper bound,
while the dashed and dotted lines represent
the separate output and input bounds,
respectively.
}
\label{fig:sharp-frontier}
\end{figure}

The resulting tight upper bound is illustrated in
Figure~\ref{fig:sharp-frontier}, together with
the separate input
and output bounds from
Section~\ref{sec:separate-bounds}.

\subsection{Exact one-step stationarity-transfer factor}
\label{subsec:one-step-factor}

The preceding results establish tight bounds relating the
normalized stationarity measures $X$ and $Z$. We now use this
relationship to determine the largest possible ratio $\frac{Z}{X}$, which
quantifies the worst-case stationarity increase across a single
inexact proximal-gradient step.

\begin{definition}[One-step stationarity-transfer factor]
\label{def:C1}
For $\eta\geq0$, define
$$
C_1^{\star}(\alpha,L,\rho,\eta)
=
\sup
\frac{\mathcal{D}_{\zeta}(\z)}
     {\mathcal{D}_{\zeta}(\x)},
$$
where the supremum is taken over all dimensions $n\geq2$,
functions satisfying the standing assumptions, and points
$\x,\z$ and residuals $\e$ satisfying
$\G\neq\mathbf{0}$,
$\e\in\mathcal{R}_{\alpha}(\x,\z)$,
$\|\e\|_2\leq\eta\|\G\|_2$, and
$\mathcal{D}_{\zeta}(\x)>0$.
\end{definition}

\begin{theorem}[Exact one-step stationarity-transfer factor]
\label{thm:C1-exact}
Suppose that $0\leq\eta<a_{\alpha}$. Then
\begin{align}\label{eq:C1-max}
C_1^{\star}(\alpha,L,\rho,\eta)
=
\max_{s\in[-\eta,\eta]}
\frac{
\theta_{\alpha}
+
\sqrt{c_{\alpha}^2+2c_{\alpha}s+\eta^2}
}{
a_{\alpha}+s
}.
\end{align}

Define
$
J_{\alpha}(\eta)
=
c_{\alpha}(a_{\alpha}-\eta)
-
|c_{\alpha}-\eta|
\bigl(\theta_{\alpha}+|c_{\alpha}-\eta|\bigr).
$
Then
\begin{align}\label{eq:C1-closed}
C_1^{\star}
=
\begin{cases}
\displaystyle
\frac{\theta_{\alpha}+|c_{\alpha}-\eta|}
     {a_{\alpha}-\eta},
& J_{\alpha}(\eta)\leq0,\\[4mm]
\displaystyle
\frac{c_{\alpha}}
{
\sqrt{1-\alpha\rho(2-\alpha L)-\eta^2}
-\theta_{\alpha}
},
& J_{\alpha}(\eta)>0.
\end{cases}
\end{align}
If $J_{\alpha}(\eta)\leq0$, the maximum in
\eqref{eq:C1-max} occurs at $s=-\eta$.
Otherwise, the maximizer is unique and belongs to
$(-\eta,\eta)$. The factor $C_1^{\star}$ is tight in every fixed dimension $n\geq2$.
\end{theorem}

\begin{proof}
To prove \eqref{eq:C1-max}, define
$$
\Phi(s)
=
\frac{\theta_{\alpha}
+\sqrt{c_{\alpha}^2+2c_{\alpha}s+\eta^2}}
{a_{\alpha}+s},
\qquad s\in[-\eta,\eta].
$$
By Definition~\ref{def:C1}, $C_1^\star$ is the
supremum of $\frac{\mathcal{D}_{\zeta}(\z)}{\mathcal{D}_{\zeta}(\x)}$ over
all admissible problems and steps.
Since $|s|\leq q\leq\eta$ and the numerator in
Corollary~\ref{cor:residual-transfer} is
nondecreasing in $q$, we have
$\frac{\mathcal{D}_{\zeta}(\z)}{\mathcal{D}_{\zeta}(\x)}
\leq\Phi(s)$.
Taking the supremum over all admissible problems
and steps yields
$$
C_1^\star
\leq\max_{s\in[-\eta,\eta]}\Phi(s).
$$

Conversely, fix any $s\in[-\eta,\eta]$.
For every $\delta>0$,
Proposition~\ref{prop:frontier-construction}
constructs an admissible problem and step satisfying
$$
\frac{\mathcal{D}_{\zeta}(\z)}
     {\mathcal{D}_{\zeta}(\x)}
=
\frac{a_{\alpha}+s}{a_{\alpha}+s+\delta}\Phi(s).
$$
Since each constructed example is admissible, we have
$
\frac{\mathcal{D}_{\zeta}(\z)}
     {\mathcal{D}_{\zeta}(\x)}
\leq C_1^\star.
$
Letting $\delta\downarrow0$ and using
$a_{\alpha}+s>0$, we obtain
$C_1^\star\geq\Phi(s)$.
Since $s\in[-\eta,\eta]$ is arbitrary, it follows that
$$
C_1^\star
\geq\max_{s\in[-\eta,\eta]}\Phi(s).
$$
Combining both inequalities proves \eqref{eq:C1-max}.

We first consider $\eta=0$. In this case,
the interval in \eqref{eq:C1-max} reduces to
$s=0$, and
$
C_1^\star
=
\frac{\theta_{\alpha}+|c_{\alpha}|}{a_{\alpha}}.
$
Since $J_{\alpha}(0)\leq0$, this gives the first case of \eqref{eq:C1-closed}.

For $\eta>0$, we distinguish two cases.
If $\alpha L\geq2$, then $c_{\alpha}\leq0$.
The numerator of $\Phi(s)$ is nonincreasing in $s$,
while its denominator is positive and increasing.
Thus, the maximum occurs at $s=-\eta$.
Since $c_{\alpha}\leq0$ and $a_{\alpha}-\eta>0$,
we have $J_{\alpha}(\eta)\leq0$.
Thus, the first case of \eqref{eq:C1-closed} follows.

Now suppose that $\alpha L<2$, so $c_{\alpha}>0$.
Define
$$
R(s)=\sqrt{c_{\alpha}^2+2c_{\alpha}s+\eta^2},
\qquad
H(s)=c_{\alpha}(a_{\alpha}+s)
-R(s)\bigl(\theta_{\alpha}+R(s)\bigr).
$$
Since $\eta>0$, we have $R(s)>0$ for
$s\in(-\eta,\eta)$. Differentiating gives
\begin{align*}
\Phi'(s)
&=\frac{H(s)}{R(s)(a_{\alpha}+s)^2},\\
H'(s)
&=-c_{\alpha}
\left(1+\frac{\theta_{\alpha}}{R(s)}\right)<0.
\end{align*}
Moreover,
$
H(-\eta)=J_{\alpha}(\eta)
$
and
$
H(\eta)=-\alpha\rho c_{\alpha}-\eta(1+\eta)<0.
$

If $J_{\alpha}(\eta)\leq0$, then
$H(-\eta)\leq0$. Since $H(s)$ is strictly
decreasing, we have $H(s)<0$ for
$s\in(-\eta,\eta)$. Hence $\Phi(s)$ is
decreasing, and its maximum occurs at
$s=-\eta$, giving the first case of
\eqref{eq:C1-closed}.

If $J_{\alpha}(\eta)>0$, then
$H(-\eta)>0$ and $H(\eta)<0$.
By continuity and strict monotonicity,
there is a unique $s^\star\in(-\eta,\eta)$
such that $H(s^\star)=0$.
Since $\Phi'(s)$ has the same sign as $H(s)$,
the function $\Phi(s)$ increases before
$s^\star$ and decreases afterward.
Thus, $s^\star$ is the unique maximizer.

At the interior maximizer $s^\star$, the identity
$H(s^\star)=0$ gives
$
\Phi(s^\star)=\frac{c_{\alpha}}{R(s^\star)}.
$
Moreover, by the definition of $R(s)$, we have
$
2c_{\alpha}s^\star
=
R(s^\star)^2-c_{\alpha}^2-\eta^2.
$
Substituting this into $H(s^\star)=0$ and
rearranging yields
$$
\begin{aligned}
\bigl(R(s^\star)+\theta_{\alpha}\bigr)^2
&=
2c_{\alpha}a_{\alpha}
+\theta_{\alpha}^2-c_{\alpha}^2-\eta^2\\
&=
1-\alpha\rho(2-\alpha L)-\eta^2.
\end{aligned}
$$
Since $R(s^\star)>0$, taking the positive
square root and substituting into
$\Phi(s^\star)=\frac{c_{\alpha}}{R(s^\star)}$
yields the second case of \eqref{eq:C1-closed}.
\end{proof}

Recall that $\kappa_{\alpha}=\max\{1,\alpha L-1\}$
was introduced in Section~\ref{sec:separate-bounds}. The following corollary compares $C_1^\star$ with
the factor obtained by combining the separate bounds
on $\mathcal{D}_{\zeta}(\x)$ and
$\mathcal{D}_{\zeta}(\z)$ from
Section~\ref{sec:separate-bounds}.

\begin{corollary}[Comparison with the separate bounds]
\label{cor:separate-comparison}
Suppose that $0<\eta<a_{\alpha}$. Then
$$
C_1^{\star}(\alpha,L,\rho,\eta)
\leq
\frac{\kappa_{\alpha}+\eta}{a_{\alpha}-\eta}.
$$
The inequality is strict if and only if $\alpha L<2$.
\end{corollary}

\begin{proof}
Suppose that $\alpha L<2$, so $c_{\alpha}>0$
and $\kappa_{\alpha}=1$.
In \eqref{eq:C1-max}, the numerator attains its
maximum $\kappa_{\alpha}+\eta$ only at $s=\eta$,
while the denominator attains its minimum
$a_{\alpha}-\eta$ only at $s=-\eta$.
Since $\eta>0$, these two extrema cannot occur
simultaneously. As the ratio is continuous on
$[-\eta,\eta]$, its maximum is strictly smaller
than $\frac{\kappa_{\alpha}+\eta}{a_{\alpha}-\eta}$.

If $\alpha L\geq2$, then $c_{\alpha}\leq0$
and $\kappa_{\alpha}=\alpha L-1$.
The numerator in \eqref{eq:C1-max} is
nonincreasing in $s$, so its maximum
$\kappa_{\alpha}+\eta$ and the minimum of
the denominator both occur at $s=-\eta$.
Hence
$
C_1^\star
=
\frac{\kappa_{\alpha}+\eta}{a_{\alpha}-\eta},
$
which proves equality.
\end{proof}

We next determine the exact condition on the prescribed
relative tolerance $\eta$ for which $C_1^\star$ is finite.

\begin{theorem}[Finite-transfer threshold]
\label{thm:transfer-threshold}
For $\eta\geq0$, we have
$$
C_1^{\star}(\alpha,L,\rho,\eta)<+\infty
\quad\Longleftrightarrow\quad
\eta<1-\alpha\rho.
$$
The same threshold holds in every fixed dimension $n\geq2$.
\end{theorem}

\begin{proof}
If $\eta<a_{\alpha}$, finiteness follows from
Theorem~\ref{thm:C1-exact}.
Now suppose that $\eta\geq a_{\alpha}$.

First, consider $a_{\alpha}>0$.
For $0<\varepsilon<a_{\alpha}$, apply
Proposition~\ref{prop:frontier-construction}
with $s=-a_{\alpha}+\varepsilon$ and
$\delta=\varepsilon$.
Since $\eta\geq a_{\alpha}$, we have
$$
\begin{aligned}
\eta^2-s^2
\geq a_{\alpha}^2-(a_{\alpha}-\varepsilon)^2
=2a_{\alpha}\varepsilon-\varepsilon^2
\geq a_{\alpha}\varepsilon.
\end{aligned}
$$
Using
$
c_{\alpha}^2+2c_{\alpha}s+\eta^2
=(c_{\alpha}+s)^2+\eta^2-s^2
$
and $\theta_{\alpha}\geq0$, we obtain
$$
\frac{\mathcal{D}_{\zeta}(\z)}
     {\mathcal{D}_{\zeta}(\x)}
=
\frac{\theta_{\alpha}
+\sqrt{c_{\alpha}^2+2c_{\alpha}s+\eta^2}}
{2\varepsilon}
\geq
\frac{\sqrt{a_{\alpha}}}{2\sqrt{\varepsilon}}
\longrightarrow+\infty
$$
as $\varepsilon\downarrow0$.

Next, consider $a_{\alpha}\leq0$.
Applying Proposition~\ref{prop:frontier-construction}
with $s=0$ and $\delta=-a_{\alpha}+\varepsilon$
gives
$$
\frac{\mathcal{D}_{\zeta}(\z)}
     {\mathcal{D}_{\zeta}(\x)}
=
\frac{\theta_{\alpha}
+\sqrt{c_{\alpha}^2+\eta^2}}{\varepsilon}
\geq
\frac{\theta_{\alpha}+|c_{\alpha}|}{\varepsilon}
\geq\frac{1}{\varepsilon}.
$$
Hence the ratio again tends to $+\infty$
as $\varepsilon\downarrow0$.
These constructions apply in every fixed dimension
$n\geq2$ and show that $C_1^\star=+\infty$
whenever $\eta\geq a_{\alpha}$.
\end{proof}

\begin{figure}[htbp]
\centering
\includegraphics[width=0.65\linewidth]
{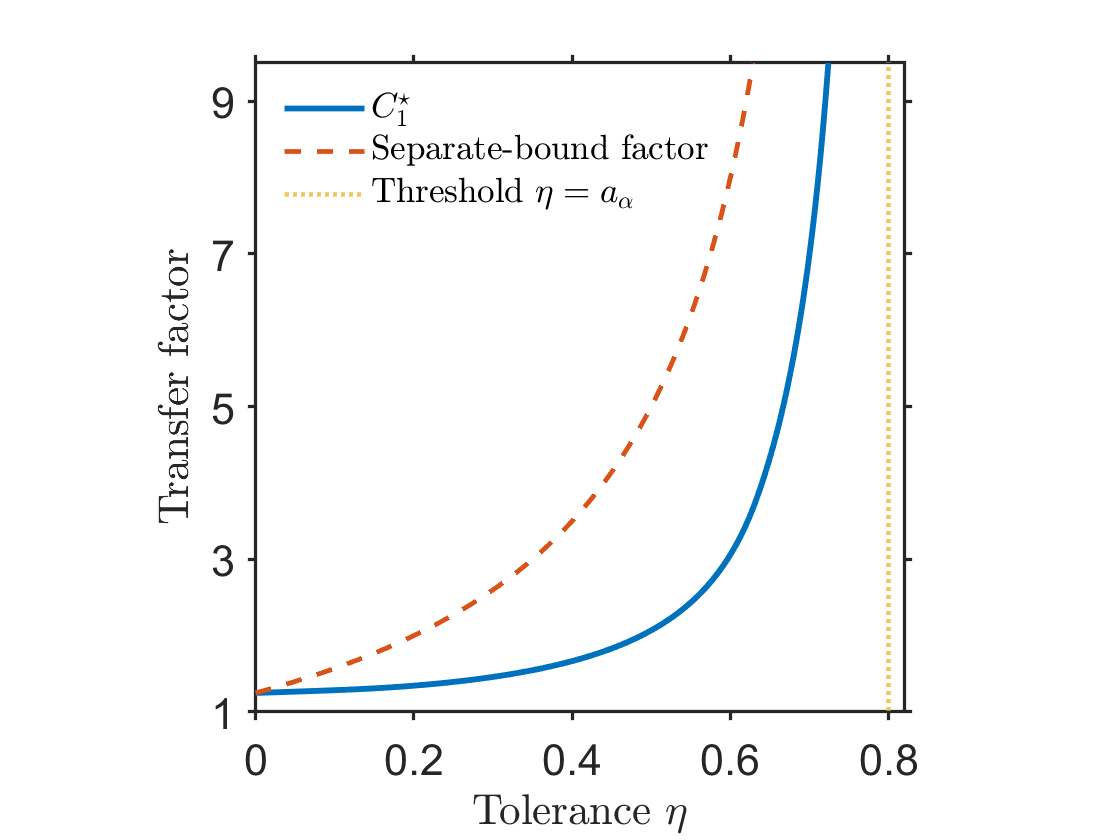}
\caption{
Comparison of the exact one-step factor
$C_1^\star$ and the separate-bound factor
$\frac{1+\eta}{0.8-\eta}$ for
$\alpha L=1$ and $\alpha\rho=0.2$.
The exact factor is strictly smaller for
$0<\eta<0.8$, while both factors diverge
as $\eta\uparrow0.8$.
}
\label{fig:one step factor}
\end{figure}

Figure~\ref{fig:one step factor} compares
$C_1^\star$ with the factor obtained from the
separate bounds, illustrating the results of
Corollary~\ref{cor:separate-comparison} and
Theorem~\ref{thm:transfer-threshold}.

\section{Repeatability and multi-step stationarity transfer}
\label{sec:repeatability}

Theorem~\ref{thm:C1-exact} establishes the exact
worst-case one-step stationarity-transfer factor
$C_1^\star$. We now investigate whether this
factor can be repeated over consecutive inexact
proximal-gradient steps of the same objective
function. Throughout this section, assume
$0\leq\eta<a_\alpha$.

\begin{definition}[Multi-step stationarity-transfer factor]
\label{def:CN}
For fixed $\alpha,L,\rho,\eta$ and
an integer $N\geq1$, define
$$
C_N^\star(\alpha,L,\rho,\eta)
=
\sup
\frac{\mathcal{D}_{\zeta}(\x^N)}
     {\mathcal{D}_{\zeta}(\x^0)},
$$
where the supremum is taken over all dimensions
$n\geq2$, all admissible pairs $(f,g)$ satisfying
the standing assumptions, and all sequences
$\{\x^k\}_{k=0}^N$ generated using the same pair
$(f,g)$ at every iteration, with
$$
\x^{k+1}=\x^k-\alpha\G_k,\quad
\e_k\in\mathcal{R}_{\alpha}(\x^k,\x^{k+1}),\quad
\|\e_k\|_2\leq\eta\|\G_k\|_2
$$
for $k=0,1,\ldots,N-1$, and
$\mathcal D_\zeta(\x^0)>0$.
For $N=1$, this definition reduces to
Definition~\ref{def:C1}.
\end{definition}

Since $C_N^\star$ bounds stationarity transfer
over any $N$ consecutive admissible steps,
splitting a sequence into two segments
yields the following property.

\begin{proposition}[Submultiplicativity]
\label{prop:submultiplicativity}
For all integers $M,N\geq1$,
$$
C_{M+N}^\star\leq C_M^\star C_N^\star.
$$
Consequently, the asymptotic factor
\begin{align}\label{eq:asymptotic-factor}
\Gamma^\star
:=
\lim_{N\to\infty}(C_N^\star)^{1/N}
=
\inf_{N\geq1}(C_N^\star)^{1/N}
\end{align}
is well-defined.
\end{proposition}

\begin{proof}
Consider any admissible $(M+N)$-step sequence
$\{\x^k\}_{k=0}^{M+N}$.
Suppose first that $\mathcal D_\zeta(\x^M)>0$.
The first $M$ steps form an admissible sequence.
For the remaining $N$ steps, define
$\y^j=\x^{M+j}$ for $j=0,1,\ldots,N$.
The reindexed sequence is also admissible
with the same objective and parameters.
Applying Definition~\ref{def:CN} to both
sequences gives
$$
\begin{aligned}
\frac{\mathcal D_\zeta(\x^{M+N})}
     {\mathcal D_\zeta(\x^0)}
=
\frac{\mathcal D_\zeta(\y^N)}
     {\mathcal D_\zeta(\y^0)}
\frac{\mathcal D_\zeta(\x^M)}
     {\mathcal D_\zeta(\x^0)}
\leq C_N^\star C_M^\star.
\end{aligned}
$$
If $\mathcal D_\zeta(\x^M)=0$, \eqref{eq:joint-input}
in Theorem~\ref{thm:joint-bounds} and
$a_\alpha-\eta>0$ imply $\G_M=\mathbf{0}$.
Hence $\x^{M+1}=\x^M$, and the same argument
applies to all subsequent steps.
Thus $\mathcal D_\zeta(\x^{M+N})=0$,
so the inequality still holds.
Taking the supremum proves submultiplicativity.

By submultiplicativity and Theorem~\ref{thm:transfer-threshold},
$C_N^\star\leq(C_1^\star)^N<\infty$.
Moreover, $C_N^\star\geq0$ by definition.
If $C_{N_0}^\star=0$ for some $N_0$,
submultiplicativity implies $C_N^\star=0$
for all $N>N_0$.
Otherwise, $\log C_N^\star$ is subadditive,
and Fekete's lemma gives
$$
\lim_{N\to\infty}(C_N^\star)^{1/N}
=
\inf_{N\geq1}(C_N^\star)^{1/N}.
$$
Thus, \eqref{eq:asymptotic-factor} holds
in either case.
\end{proof}

Proposition~\ref{prop:submultiplicativity}
implies $C_N^\star\leq(C_1^\star)^N$. We first identify a regime in which
this upper bound is attained for every $N\geq1$.

\begin{theorem}[Exact endpoint repeatability]
\label{thm:endpoint-repeatability}
Suppose that $J_{\alpha}(\eta)\leq0$.
Then, for every integer $N\geq1$,
\begin{align}\label{eq:endpoint-repeatability}
C_N^\star=(C_1^\star)^N.
\end{align}
\end{theorem}

\begin{proof}
By Proposition~\ref{prop:submultiplicativity},
$C_N^\star\leq(C_1^\star)^N$.
Since $J_\alpha(\eta)\leq0$,
Theorem~\ref{thm:C1-exact} shows that the
maximum in \eqref{eq:C1-max} is attained
at $s=-\eta$, with
$
C_1^\star
=
\frac{\theta_\alpha+|c_\alpha-\eta|}
     {a_\alpha-\eta}.
$
It therefore suffices to construct an
admissible sequence for which this factor
is attained at every step.

In both cases, set $g(t)=-\frac{\rho}{2}t^2$
and define $\omega_k:=\zeta'(x^k)$.
To attain the endpoint value $s=-\eta$,
we require $e_k=-\eta G_k$.
Since $x^{k+1}=x^k-\alpha G_k$,
the residual identity gives
$$
e_k=g'(x^{k+1})+f'(x^k)-G_k
=\omega_k-a_\alpha G_k.
$$
Combining this identity with $e_k=-\eta G_k$
yields
$G_k=\frac{\omega_k}{a_\alpha-\eta}$,
which is well-defined since $a_\alpha-\eta>0$.
Thus, the constructed residual vectors satisfy
$e_k\in\mathcal{R}_{\alpha}(x^k,x^{k+1})$,
and the relative residual condition \eqref{eq:relative-residual} holds with
equality at every step.

If $c_\alpha\geq\eta$, take $f(t)=\ell t$
with $\ell\neq0$ and choose $x^0$ so that
$\omega_0\neq0$.
Since $\omega_k=\ell-\rho x^k$,
we have
$$
\begin{aligned}
\omega_{k+1}
=\omega_k+\alpha\rho G_k
=\frac{1-\eta}{a_\alpha-\eta}\,\omega_k
=C_1^\star\omega_k.
\end{aligned}
$$

If $c_\alpha<\eta$, take $f(t)=\frac{Lt^2}{2}$
and choose $x^0$ with $\omega_0\neq0$.
Since $\omega_k=(L-\rho)x^k$,
we obtain
$$
\begin{aligned}
\omega_{k+1}
=\omega_k-\alpha(L-\rho)G_k
=\frac{1-\eta-\alpha L}{a_\alpha-\eta}\,
\omega_k
=-C_1^\star\omega_k.
\end{aligned}
$$

In both cases, $\mathcal D_\zeta(x^k)=|\omega_k|$,
so the preceding identities imply
$\mathcal D_\zeta(x^{k+1})
=C_1^\star\mathcal D_\zeta(x^k)$
for every $k$. Consequently,
$$
\frac{\mathcal D_\zeta(x^N)}
     {\mathcal D_\zeta(x^0)}
=(C_1^\star)^N.
$$
The scalar construction extends to every
dimension $n\geq2$ by letting $f$ and $g$
depend only on the first coordinate and
performing all updates along that coordinate.
By Definition~\ref{def:CN}, the construction
yields $C_N^\star\geq(C_1^\star)^N$, completing the proof.
\end{proof}

The endpoint construction uses
$\e_k=-\eta\G_k$, so the residual is
parallel to the step direction.
By contrast, an interior one-step maximizer
requires a nonzero residual component
orthogonal to the step direction.
The following lemma establishes necessary
conditions for steps whose stationarity-transfer
ratios approach $C_1^\star$ in this regime.

\begin{lemma}[Structure of interior near-maximizers]
\label{lem:interior-rigidity}
Suppose that $J_\alpha(\eta)>0$, and let
$s^\star\in(-\eta,\eta)$ be the unique maximizer
in \eqref{eq:C1-max}. Set
$$
r^\star=\sqrt{\eta^2-(s^\star)^2}>0,
\qquad
Q^\star=\sqrt{c_\alpha^2+2c_\alpha s^\star+\eta^2}>0.
$$
Consider a sequence of admissible one-step
instances $\x_j\to\z_j$, indexed by $j\geq1$, with objectives
$\zeta_j=f_j+g_j$ and possibly different
dimensions, such that
$
\frac{\mathcal D_{\zeta_j}(\z_j)}
     {\mathcal D_{\zeta_j}(\x_j)}
\to C_1^\star.
$
For each instance, define
$t_j=\|\G_j\|_2$ and $\bu_j=\frac{\G_j}{t_j}$.

Let $\w_{x,j}$ and $\w_{z,j}$ denote the
minimum-norm elements of
$S_{\zeta_j}(\x_j)$ and $S_{\zeta_j}(\z_j)$,
respectively.
By the residual condition,
$
\widehat{\w}_{z,j}
:=
\G_j+\e_j+\nabla f_j(\z_j)
-\nabla f_j(\x_j)
\in S_{\zeta_j}(\z_j).
$
Then there exist unit vectors
$\bv_j\perp\bu_j$ such that
$$
\begin{aligned}
\frac{\e_j}{t_j}-s^\star\bu_j-r^\star\bv_j
&\to\mathbf{0},
&\quad
\frac{\w_{x,j}}{t_j}-(a_\alpha+s^\star)\bu_j
&\to\mathbf{0},\\
\frac{\widehat{\w}_{z,j}-\w_{z,j}}{t_j}
&\to\mathbf{0},&\quad
\frac{\w_{z,j}}{\|\w_{z,j}\|_2}
-\frac{(c_\alpha+s^\star)\bu_j+r^\star\bv_j}{Q^\star}
&\to\mathbf{0}.
\end{aligned}
$$
All vector convergence statements mean
convergence to zero in Euclidean norm as $j\to\infty$.
\end{lemma}

\begin{proof}
For the $j$th instance, let $q_j,s_j$ be the
residual quantities, and define
$X_j=\frac{\mathcal D_{\zeta_j}(\x_j)}{t_j}$,
$Z_j=\frac{\mathcal D_{\zeta_j}(\z_j)}{t_j}$,
and $R_j=\frac{Z_j}{X_j}$.
By Corollary~\ref{cor:residual-transfer},
$$
R_j
\leq
\frac{\theta_\alpha+
\sqrt{c_\alpha^2+2c_\alpha s_j+q_j^2}}
{a_\alpha+s_j}
\leq C_1^\star.
$$
The residual constraints $|s_j|\leq q_j\leq\eta$
ensure that the pairs $(q_j,s_j)$ lie in a
compact set.
The fraction in the preceding inequality is continuous on this
set and has the unique maximizer
$(q,s)=(\eta,s^\star)$.
Since $R_j\to C_1^\star$, continuity,
compactness, and uniqueness imply
$q_j\to\eta$ and $s_j\to s^\star$.

Set $r_j=\sqrt{q_j^2-s_j^2}$, so that
$r_j\to r^\star>0$.
For all sufficiently large $j$, we can
decompose the scaled residual as
$\frac{\e_j}{t_j}=s_j\bu_j+r_j\bv_j$, where $\bv_j$ is a unit vector orthogonal
to $\bu_j$.
Since $s_j\to s^\star$ and
$r_j\to r^\star$, the first limit follows.

To establish the remaining limits, the joint bounds \eqref{eq:joint-input}
and \eqref{eq:joint-output} give
$$
a_\alpha+s_j
\leq X_j
=\frac{Z_j}{R_j}
\leq
\frac{\theta_\alpha+
\sqrt{c_\alpha^2+2c_\alpha s_j+q_j^2}}
{R_j}.
$$
Using $s_j\to s^\star$, $q_j\to\eta$,
and $R_j\to C_1^\star$, the lower bound
converges to $a_\alpha+s^\star$, while
the upper bound converges to
$\frac{\theta_\alpha+Q^\star}{C_1^\star}
=a_\alpha+s^\star$.
Hence, the squeeze theorem gives
$X_j\to a_\alpha+s^\star$.
The identity $Z_j=R_jX_j$ then yields
$Z_j\to\theta_\alpha+Q^\star$.

For the second limit, the proof of
Theorem~\ref{thm:joint-bounds} and
the Cauchy-Schwarz inequality give
$$
a_\alpha+s_j
\leq
\left\langle
\frac{\w_{x,j}}{t_j},\bu_j
\right\rangle
\leq
\frac{\|\w_{x,j}\|_2}{t_j}
=X_j.
$$
As $s_j\to s^\star$ and
$X_j\to a_\alpha+s^\star$, the squeeze
theorem gives
$
\left\langle
\frac{\w_{x,j}}{t_j},\bu_j
\right\rangle
\to a_\alpha+s^\star.
$
Combining this with
$\frac{\|\w_{x,j}\|_2}{t_j}\to a_\alpha+s^\star$
and $\|\bu_j\|_2=1$ yields
$
\left\|
\frac{\w_{x,j}}{t_j}
-(a_\alpha+s^\star)\bu_j
\right\|_2\to0,
$
which proves the second limit.

Next, define
$$
\mathbf b_j=c_\alpha\bu_j+\frac{\e_j}{t_j},
\qquad
\boldsymbol\beta_j
=
-\frac{\nabla f_j(\z_j)-\nabla f_j(\x_j)}{t_j}
-\theta_\alpha\bu_j.
$$
By Lemma~\ref{lem:coco-ball},
$\|\boldsymbol\beta_j\|_2\leq\theta_\alpha$.
The definition of $\widehat{\w}_{z,j}$ gives
$\frac{\widehat{\w}_{z,j}}{t_j}
=\mathbf b_j-\boldsymbol\beta_j$,
while the first limit implies
$\|\mathbf b_j\|_2\to Q^\star>0$.

Since $\widehat{\w}_{z,j}\in
S_{\zeta_j}(\z_j)$ and $\w_{z,j}$
has minimum norm,
$$
Z_j
\leq\frac{\|\widehat{\w}_{z,j}\|_2}{t_j}
=\|\mathbf b_j-\boldsymbol\beta_j\|_2
\leq\|\mathbf b_j\|_2+\theta_\alpha.
$$
The limits $Z_j\to\theta_\alpha+Q^\star$
and $\|\mathbf b_j\|_2\to Q^\star$
then yield
$\frac{\|\widehat{\w}_{z,j}\|_2}{t_j}
\to\theta_\alpha+Q^\star$
by the squeeze theorem.
Consequently,
$\|\mathbf b_j-\boldsymbol\beta_j\|_2
\to Q^\star+\theta_\alpha$.
Together with
$\|\boldsymbol\beta_j\|_2\leq\theta_\alpha$
and $\|\mathbf b_j\|_2\to Q^\star>0$,
the equality condition for the triangle
inequality gives
$
\frac{\widehat{\w}_{z,j}}
{\|\widehat{\w}_{z,j}\|_2}
-
\frac{\mathbf b_j}{\|\mathbf b_j\|_2}
\to\mathbf{0}.
$

The minimum-norm property of $\w_{z,j}$
and the convexity of $S_{\zeta_j}(\z_j)$
yield
$$
0\leq
\frac{\|\widehat{\w}_{z,j}-\w_{z,j}\|_2^2}{t_j^2}
\leq
\frac{\|\widehat{\w}_{z,j}\|_2^2
-\|\w_{z,j}\|_2^2}{t_j^2}
\to0.
$$
This proves the third limit.
Furthermore, both
$\frac{\|\widehat{\w}_{z,j}\|_2}{t_j}$ and
$\frac{\|\w_{z,j}\|_2}{t_j}$ converge to
$\theta_\alpha+Q^\star>0$, so
$
\frac{\widehat{\w}_{z,j}}
{\|\widehat{\w}_{z,j}\|_2}
-
\frac{\w_{z,j}}
{\|\w_{z,j}\|_2}
\to\mathbf{0}.
$

The first limit also gives
$
\mathbf b_j-
\bigl((c_\alpha+s^\star)\bu_j
+r^\star\bv_j\bigr)\to\mathbf{0}.
$
Together with $\|\mathbf b_j\|_2\to Q^\star$,
this establishes the fourth limit.
\end{proof}

\begin{theorem}[Strict interior nonrepeatability]
\label{thm:interior-nonrepeatability}
Suppose that $J_{\alpha}(\eta)>0$.
Then
\begin{align}\label{eq:strict-two-step}
C_2^\star<(C_1^\star)^2.
\end{align}
Thus, the exact one-step factor $C_1^\star$
cannot be attained or approached arbitrarily
closely over two consecutive inexact
proximal-gradient steps of the same objective.
\end{theorem}

\begin{proof}
By Proposition \ref{prop:submultiplicativity},
$C_2^\star\leq(C_1^\star)^2$.
To prove that this inequality is strict,
suppose that $C_2^\star=(C_1^\star)^2$.
Then there exists a sequence of admissible
two-step instances whose stationarity-transfer
ratios converge to $(C_1^\star)^2$.
Since each one-step ratio is at most
$C_1^\star$, both one-step ratios must
converge to $C_1^\star$.
We can therefore apply
Lemma~\ref{lem:interior-rigidity}
to each step.

For the $j$th two-step instance, let
$k=0,1$ index the two consecutive steps,
and define
$t_{k,j}=\|\G_{k,j}\|_2$ and
$\bu_{k,j}=\frac{\G_{k,j}}{t_{k,j}}$.
Let $\bv_{k,j}\perp\bu_{k,j}$ denote
the unit vectors given by
Lemma~\ref{lem:interior-rigidity}.
The point $\x_j^1$ is the output of the
first step and the input of the second.
Thus, both steps share the same minimum-norm
element of $S_{\zeta_j}(\x_j^1)$.
Applying the fourth limit of
Lemma~\ref{lem:interior-rigidity} to the
first step and the second limit to the
second step at the shared point $\x_j^1$
gives
$$
\bu_{1,j}
-
\frac{(c_\alpha+s^\star)\bu_{0,j}
+r^\star\bv_{0,j}}{Q^\star}
\to\mathbf{0}.
$$
Since $\bv_{0,j}\perp\bu_{0,j}$, we obtain
$\langle\bv_{0,j},\bu_{1,j}\rangle
\to \frac{r^\star}{Q^\star}>0$.

Let $\widehat g_j=g_j+\frac{\rho\|\cdot\|_2^2}{2}$.
For $k=0,1$, let $\w_{k,j}$ denote the
minimum-norm element of
$S_{\zeta_j}(\x_j^k)$, and define
$
\mathbf p_{k,j}
:=\w_{k,j}-\nabla f_j(\x_j^k)+\rho\x_j^k \text{ and }
\widehat{\mathbf p}_{k+1,j}
:=\G_{k,j}+\e_{k,j}
-\nabla f_j(\x_j^k)+\rho\x_j^{k+1}.
$
The first vector belongs to
$\partial\widehat g_j(\x_j^k)$,
while the second belongs to
$\partial\widehat g_j(\x_j^{k+1})$
by the residual condition.
At the final point, set
$\mathbf p_{2,j}:=\widehat{\mathbf p}_{2,j}$.
Using the definitions above and
$\x_j^{k+1}-\x_j^k=-\alpha\G_{k,j}$,
we obtain
$$
\widehat{\mathbf p}_{k+1,j}-\mathbf p_{k,j}
=a_\alpha\G_{k,j}+\e_{k,j}-\w_{k,j}.
$$

For the first step, the third limit of
Lemma~\ref{lem:interior-rigidity} implies
$\frac{\widehat{\mathbf p}_{1,j}
-\mathbf p_{1,j}}{t_{0,j}}\to\mathbf{0}$,
since this difference equals the difference
between the corresponding stationarity
vectors at $\x_j^1$.
Applying the first two limits of the lemma
to both steps then yields
$$
\begin{aligned}
\frac{\mathbf p_{1,j}-\mathbf p_{0,j}}{t_{0,j}}
-r^\star\bv_{0,j}
\to\mathbf{0},\quad
\frac{\mathbf p_{2,j}-\mathbf p_{1,j}}{t_{1,j}}
-r^\star\bv_{1,j}
\to\mathbf{0}.
\end{aligned}
$$
At the shared point $\x_j^1$, applying the
limits for $Z_j$ and $X_j$ established in
the proof of Lemma~\ref{lem:interior-rigidity}
to the first and second steps, respectively,
gives
$$
\frac{t_{1,j}}{t_{0,j}}
=
\frac{\mathcal D_{\zeta_j}(\x_j^1)/t_{0,j}}
     {\mathcal D_{\zeta_j}(\x_j^1)/t_{1,j}}
\to
\frac{\theta_\alpha+Q^\star}
     {a_\alpha+s^\star}
=C_1^\star.
$$

The convexity of $\widehat g_j$ implies that
$\partial\widehat g_j$ is cyclically monotone \cite{rockafellar1998variational}.
Applying this property to the cycle
$\x_j^0\to\x_j^2\to\x_j^1\to\x_j^0$
gives
$$
\begin{aligned}
0\leq{}
\langle\mathbf p_{0,j},\x_j^0-\x_j^2\rangle
+\langle\mathbf p_{2,j},\x_j^2-\x_j^1\rangle+\langle\mathbf p_{1,j},\x_j^1-\x_j^0\rangle.
\end{aligned}
$$
Using
$\x_j^{k+1}-\x_j^k=-\alpha t_{k,j}\bu_{k,j}$
and dividing by $\alpha t_{0,j}t_{1,j}>0$,
we obtain
$$
\begin{aligned}
0\leq{}
-\frac{t_{0,j}}{t_{1,j}}
\left\langle
\frac{\mathbf p_{1,j}-\mathbf p_{0,j}}{t_{0,j}},
\bu_{0,j}
\right\rangle-
\left\langle
\frac{\mathbf p_{2,j}-\mathbf p_{0,j}}{t_{0,j}},
\bu_{1,j}
\right\rangle.
\end{aligned}
$$

The convergence relation for
$\frac{\mathbf p_{1,j}-\mathbf p_{0,j}}{t_{0,j}}$,
together with $\bv_{0,j}\perp\bu_{0,j}$ and
$\frac{t_{0,j}}{t_{1,j}}\to\frac{1}{C_1^\star}$,
shows that the first term tends to zero.
For the second term, we have
$$
\begin{aligned}
\left\langle
\frac{\mathbf p_{2,j}-\mathbf p_{0,j}}{t_{0,j}},
\bu_{1,j}
\right\rangle
=
\left\langle
\frac{\mathbf p_{1,j}-\mathbf p_{0,j}}{t_{0,j}},
\bu_{1,j}
\right\rangle
+
\frac{t_{1,j}}{t_{0,j}}
\left\langle
\frac{\mathbf p_{2,j}-\mathbf p_{1,j}}{t_{1,j}},
\bu_{1,j}
\right\rangle
\to \frac{(r^\star)^2}{Q^\star},
\end{aligned}
$$
where we used
$\langle\bv_{0,j},\bu_{1,j}\rangle
\to \frac{r^\star}{Q^\star}$
and $\bv_{1,j}\perp\bu_{1,j}$.
Consequently, the right-hand side of
the cyclic-monotonicity inequality converges
to $-\frac{(r^\star)^2}{Q^\star}<0$,
contradicting its nonnegativity.
This establishes \eqref{eq:strict-two-step}.
\end{proof}

\begin{corollary}[Asymptotic stationarity-transfer factor]
\label{cor:persistent-transient}
Let $\Gamma^\star$ be defined by
\eqref{eq:asymptotic-factor}, and define
the left-endpoint factor
$C_-=\frac{\theta_\alpha+|c_\alpha-\eta|}{a_\alpha-\eta}$.
If $J_\alpha(\eta)\leq0$, then
$\Gamma^\star=C_1^\star$.
If $J_\alpha(\eta)>0$, then
$$
C_-\leq\Gamma^\star
\leq\sqrt{C_2^\star}<C_1^\star.
$$
\end{corollary}

\begin{proof}
If $J_\alpha(\eta)\leq0$,
Theorem~\ref{thm:endpoint-repeatability} gives
$C_N^\star=(C_1^\star)^N$ for every $N\geq1$.
Thus $\Gamma^\star=C_1^\star$ by
\eqref{eq:asymptotic-factor}.

Now suppose that $J_\alpha(\eta)>0$.
The scalar construction in the proof of
Theorem~\ref{thm:endpoint-repeatability}
remains admissible in this regime and
attains the factor $C_-$ at every step.
Consequently, $C_N^\star\geq C_-^N$
for every $N\geq1$, so
$(C_N^\star)^{\frac{1}{N}}\geq C_-$.
By \eqref{eq:asymptotic-factor},
we obtain $\Gamma^\star\geq C_-$.
On the other hand,
Proposition~\ref{prop:submultiplicativity}
and Theorem~\ref{thm:interior-nonrepeatability}
give
$$
\Gamma^\star
\leq\sqrt{C_2^\star}<C_1^\star.
$$
Combining these bounds proves the result.
\end{proof}

\section{Tolerance selection and stationarity convergence}
\label{sec:accuracy-descent}

Section~\ref{sec:repeatability} investigates
whether the worst-case one-step
stationarity-transfer factor $C_1^\star$
can be repeated over consecutive inexact
proximal-gradient steps.
We now study a complementary aspect of
multi-step behavior: stationarity convergence
along an iterative sequence.
We first use $C_1^\star$ to determine the
largest prescribed relative tolerance that
guarantees a desired stationarity-transfer bound.
We then derive a one-step descent estimate
and establish stationarity convergence and
complexity guarantees for the resulting
inexact proximal-gradient sequence.

\begin{corollary}[Largest relative tolerance]
\label{cor:largest-tolerance}
Suppose that $a_{\alpha}>0$, and set
$
C_{\mathrm{ex}}=C_1^\star(\alpha,L,\rho,0)=\frac{\kappa_{\alpha}}{a_{\alpha}}.
$
For every $\gamma>C_{\mathrm{ex}}$, there exists a unique
$\eta_{\gamma}\in(0,a_{\alpha})$ such that
$$
C_1^\star(\alpha,L,\rho,\eta_{\gamma})=\gamma.
$$
Moreover, the bound
$
\mathcal{D}_{\zeta}(\z)\le \gamma \mathcal{D}_{\zeta}(\x)
$
holds uniformly over all admissible problems and steps satisfying
$
\|\e\|_2\le \eta \|\G\|_2
$
if and only if $0\le \eta\le \eta_{\gamma}$.

If $\alpha L\ge 2$, then
$$
\eta_{\gamma}
=
\frac{\gamma a_{\alpha}-(\alpha L-1)}{\gamma+1}.
$$

If $\alpha L<2$, define
$$
\eta_1=\frac{\gamma a_{\alpha}-1}{\gamma-1},
\
\eta_2=\frac{\gamma a_{\alpha}-(\alpha L-1)}{\gamma+1},
\
\eta_{\mathrm{int}}
=
\sqrt{
1-\alpha\rho(2-\alpha L)
-\left(\theta_{\alpha}+\frac{c_{\alpha}}{\gamma}\right)^2
}.
$$
Among the real candidates in $(0,a_{\alpha})$,
the value of $\eta_{\gamma}$ is determined by
$$
\eta_{\gamma}=
\begin{cases}
\eta_1,
& \eta_1\leq c_{\alpha},
\quad J_{\alpha}(\eta_1)\leq0,\\[2mm]
\eta_2,
& \eta_2\geq c_{\alpha},
\quad J_{\alpha}(\eta_2)\leq0,\\[2mm]
\eta_{\mathrm{int}},
& J_{\alpha}(\eta_{\mathrm{int}})>0.
\end{cases}
$$
\end{corollary}

\begin{proof}
Let $C(\eta)=C_1^\star(\alpha,L,\rho,\eta)$ for
$0\leq\eta<a_{\alpha}$.
To show that $C$ is continuous, substitute
$s=\eta t$ into \eqref{eq:C1-max} to obtain
$$
C(\eta)
=
\max_{t\in[-1,1]}
\frac{\theta_{\alpha}
+\sqrt{c_{\alpha}^2+2c_{\alpha}\eta t+\eta^2}}
{a_{\alpha}+\eta t}.
$$
Since $a_{\alpha}+\eta t\geq a_{\alpha}-\eta>0$,
the expression being maximized is continuous in
$(\eta,t)$. As the maximization is over the fixed
compact interval $[-1,1]$, $C(\eta)$ is continuous.

We next show that $C$ is strictly increasing.
Fix $0\leq\eta<\widetilde{\eta}<a_{\alpha}$,
and let $s_{\eta}$ be a maximizer in
\eqref{eq:C1-max} for the tolerance $\eta$.
Since $|s_{\eta}|\leq\eta<\widetilde{\eta}$,
the same $s_{\eta}$ is also admissible for
$\widetilde{\eta}$.
For this fixed $s_{\eta}$, increasing $\eta$
to $\widetilde{\eta}$ strictly increases the
numerator in \eqref{eq:C1-max}, while the
denominator remains unchanged.
Thus, the value at $s_{\eta}$ for
$\widetilde{\eta}$ exceeds $C(\eta)$.
Since $C(\widetilde{\eta})$ is the maximum
over all admissible $s$, we conclude that
$C(\widetilde{\eta})>C(\eta)$.

We also show that $C(\eta)\to+\infty$ as
$\eta\uparrow a_{\alpha}$.
Let $\eta=a_{\alpha}-\varepsilon$ and choose
$s=-a_{\alpha}+2\varepsilon$, where
$0<\varepsilon<\frac{a_{\alpha}}{3}$.
Then $s\in[-\eta,\eta]$ and
$
\eta^2-s^2
=2a_{\alpha}\varepsilon-3\varepsilon^2
\geq a_{\alpha}\varepsilon.
$
Since
$
c_{\alpha}^2+2c_{\alpha}s+\eta^2
=(c_{\alpha}+s)^2+\eta^2-s^2
$
and $\theta_{\alpha}\geq0$,
\eqref{eq:C1-max} gives
$$
C(\eta)
\geq
\frac{\sqrt{a_{\alpha}}}{2\sqrt{\varepsilon}}
\longrightarrow+\infty
\qquad\text{as }\varepsilon\downarrow0.
$$

Since $C(0)=C_{\mathrm{ex}}$ and $C$ is continuous
and strictly increasing, for every
$\gamma>C_{\mathrm{ex}}$ there exists a unique
$\eta_{\gamma}\in(0,a_{\alpha})$ satisfying
$C(\eta_{\gamma})=\gamma$.
By Theorems~\ref{thm:C1-exact} and
\ref{thm:transfer-threshold},
$\mathcal{D}_{\zeta}(\z)\leq
\gamma\mathcal{D}_{\zeta}(\x)$
holds uniformly over all admissible problems
and steps if and only if
$0\leq\eta\leq\eta_{\gamma}$.
The stated formulas follow by solving the
corresponding branches of \eqref{eq:C1-closed}
for $\eta$.
\end{proof}

Corollary~\ref{cor:largest-tolerance} determines
the largest prescribed relative tolerance $\eta$
that guarantees a given stationarity-transfer bound.
However, this guarantee does not necessarily imply
a decrease in the objective function.
We therefore examine objective descent under the
same residual condition, retaining the residual
orientation $s$ in the analysis.

\begin{proposition}[Angle-aware descent]
\label{prop:angle-descent}
Consider an inexact proximal-gradient step
from $\x$ to $\z$ with $\G\neq\mathbf{0}$
and residual $\e\in\mathcal{R}_{\alpha}(\x,\z)$.
Then
\begin{align}\label{eq:angle-descent}
\zeta(\z)
\leq\zeta(\x)
-\alpha\left(1-\frac{\alpha(L+\rho)}{2}+s\right)
\|\G\|_2^2.
\end{align}
In particular, if $\|\e\|_2\leq\eta\|\G\|_2$,
then
$$
\zeta(\z)
\leq\zeta(\x)
-\alpha\left(1-\frac{\alpha(L+\rho)}{2}-\eta\right)
\|\G\|_2^2.
$$
\end{proposition}

\begin{proof}
Since $\e\in\mathcal{R}_{\alpha}(\x,\z)$,
there exists $\bv_{\z}\in\partial\widehat{g}(\z)$
such that
$
\nabla f(\x)+\bv_{\z}-\rho\z=\G+\e.
$
By the convexity of $\widehat{g}$, the identity
$g(\x)=\widehat{g}(\x)-\frac{\rho}{2}\|\x\|_2^2$,
and the $L$-Lipschitz continuity of $\nabla f$,
we have
$$
\begin{aligned}
g(\z)
&\leq g(\x)
+\langle\bv_{\z}-\rho\z,\z-\x\rangle
+\frac{\rho}{2}\|\z-\x\|_2^2,\\
f(\z)
&\leq f(\x)
+\langle\nabla f(\x),\z-\x\rangle
+\frac{L}{2}\|\z-\x\|_2^2.
\end{aligned}
$$
Adding these inequalities and using
$\z-\x=-\alpha\G$ gives
$$
\zeta(\z)
\leq\zeta(\x)
-\alpha\langle\G+\e,\G\rangle
+\frac{\alpha^2(L+\rho)}{2}\|\G\|_2^2.
$$
Since $\langle\e,\G\rangle=s\|\G\|_2^2$,
we obtain \eqref{eq:angle-descent}.
The second inequality follows from $s\geq-\eta$.
\end{proof}

Proposition~\ref{prop:angle-descent} establishes
an objective-descent estimate for a single
inexact proximal-gradient step, while
Theorem~\ref{thm:joint-bounds} bounds the
stationarity measure at its output.
Combining these estimates over successive steps,
while retaining their common dependence on $s$,
yields the following stationarity convergence
and complexity result.

\begin{corollary}[Stationarity convergence and complexity]
\label{cor:stationarity-complexity}
Consider a sequence of inexact proximal-gradient steps
$\x^{k+1}=\x^k-\alpha\G_k$ with fixed $\alpha>0$
such that
$\e_k\in\mathcal{R}_{\alpha}(\x^k,\x^{k+1})$
and
$\|\e_k\|_2\leq\eta\|\G_k\|_2$
for all $k$.
Suppose that $\eta\geq0$ is fixed,
$\zeta_{\inf}:=\inf_{\w\in\mathbb{R}^n}\zeta(\w)>-\infty$,
and
$
d:=1-\frac{\alpha(L+\rho)}{2}-\eta>0.
$
Define
$$
K_{\alpha}(\eta)
=
\max_{s\in[-\eta,\eta]}
\frac{
\theta_{\alpha}
+\sqrt{c_{\alpha}^2+2c_{\alpha}s+\eta^2}
}{
\sqrt{d+\eta+s}
}.
$$
If the sequence is infinite, then
$$
\mathcal{D}_{\zeta}(\x^k)\to0
\qquad\text{as }k\to\infty.
$$
Moreover, for every integer $N\geq1$ for which
the first $N$ steps are defined,
\begin{align}\label{eq:stationarity-complexity}
\min_{1\leq j\leq N}\mathcal{D}_{\zeta}(\x^j)
\leq
K_{\alpha}(\eta)
\sqrt{
\frac{\zeta(\x^0)-\zeta_{\inf}}{\alpha N}
}
\leq
(1+\eta)
\sqrt{
\frac{\zeta(\x^0)-\zeta_{\inf}}{\alpha dN}
}.
\end{align}
The coefficient satisfies
$
K_{\alpha}(\eta)<\frac{1+\eta}{\sqrt{d}}
$
whenever $\eta>0$.
\end{corollary}

\begin{proof}
For each step with $\G_k\neq\mathbf{0}$,
let $q_k$ and $s_k$ be the residual quantities
defined in Section~2.
By Theorem~\ref{thm:joint-bounds}, we have
$$
\mathcal{D}_{\zeta}(\x^{k+1})
\leq
\left(
\theta_{\alpha}
+\sqrt{c_{\alpha}^2+2c_{\alpha}s_k+q_k^2}
\right)\|\G_k\|_2.
$$
Moreover, Proposition~\ref{prop:angle-descent}
shows that the objective decrease
$\zeta(\x^k)-\zeta(\x^{k+1})$ is at least
$\alpha(d+\eta+s_k)\|\G_k\|_2^2$.
Since $|s_k|\leq q_k\leq\eta$ and
$d+\eta+s_k\geq d>0$, combining these
estimates yields
$$
\mathcal{D}_{\zeta}(\x^{k+1})^2
\leq
\frac{K_{\alpha}(\eta)^2}{\alpha}
\bigl(\zeta(\x^k)-\zeta(\x^{k+1})\bigr).
$$
If $\G_k=\mathbf{0}$, then $\x^{k+1}=\x^k$.
The relative residual condition gives
$\|\e_k\|_2\leq\eta\|\G_k\|_2=0$,
so $\e_k=\mathbf{0}$.
By the definition of $\mathcal{R}_{\alpha}$ in \eqref{eq:PG-residual},
we have
$
\mathbf{0}\in\mathcal{S}_{\zeta}(\x^{k+1}),
$
which implies
$\mathcal{D}_{\zeta}(\x^{k+1})=0$
by Lemma~\ref{lem:d-gap}.
Hence the same inequality holds.

Summing over $k=0,\ldots,N-1$ gives
$$
\sum_{j=1}^{N}
\mathcal{D}_{\zeta}(\x^j)^2
\leq
\frac{K_{\alpha}(\eta)^2}{\alpha}
\bigl(\zeta(\x^0)-\zeta_{\inf}\bigr).
$$
If the sequence is infinite, the right-hand
side is finite and independent of $N$.
Therefore,
$
\sum_{j=1}^{\infty}
\mathcal{D}_{\zeta}(\x^j)^2<+\infty,
$
which implies
$\mathcal{D}_{\zeta}(\x^k)\to0$
as $k\to\infty$.

For every $N\geq1$, the summation bound gives
$$
\begin{aligned}
\min_{1\leq j\leq N}
\mathcal{D}_{\zeta}(\x^j)^2
\leq
\frac{1}{N}\sum_{j=1}^{N}
\mathcal{D}_{\zeta}(\x^j)^2
\leq
\frac{K_{\alpha}(\eta)^2}{\alpha N}
\bigl(\zeta(\x^0)-\zeta_{\inf}\bigr).
\end{aligned}
$$
Taking square roots yields the first inequality
in \eqref{eq:stationarity-complexity}.

It remains to establish the second inequality in
\eqref{eq:stationarity-complexity}.
Since $d>0$, we have $\alpha L<2$
and hence $c_{\alpha}>0$.
Thus, the numerator in the definition of
$K_{\alpha}(\eta)$ attains its maximum
$1+\eta$ at $s=\eta$, while the denominator
attains its minimum $\sqrt{d}$ at $s=-\eta$.
It follows that
$
K_{\alpha}(\eta)\leq\frac{1+\eta}{\sqrt{d}}.
$
When $\eta>0$, these extrema cannot occur
simultaneously. Since the maximum defining
$K_{\alpha}(\eta)$ is attained on the compact
interval $[-\eta,\eta]$, the inequality is strict.
\end{proof}

\section{Conclusion}
\label{sec:conclusion}

In this paper, we investigated the relationship
between the directional-stationarity measures at
the input and output points of inexact
proximal-gradient steps for composite objectives
with a smooth convex term and a possibly
nonsmooth weakly convex term.
By retaining information about both the
relative magnitude of the residual and
its alignment with the proximal-gradient
displacement, we established joint
stationarity bounds and derived a tight input-output upper bound.
This analysis yielded the exact
worst-case one-step stationarity-transfer
factor $C_1^\star$, identified the precise
relative-residual threshold for its finiteness,
and quantified the improvement over the
factor obtained from separate input and
output bounds.

We further examined whether the worst-case
one-step factor can be repeated over consecutive
inexact proximal-gradient steps of the same
objective. In the endpoint regime, we established
exact repeatability over arbitrary numbers of
steps. In contrast, the interior regime exhibits
strict two-step nonrepeatability, leading to an
asymptotic stationarity-transfer factor strictly
smaller than $C_1^\star$.
We also determined the largest relative-residual
tolerance that guarantees a prescribed
one-step stationarity-transfer bound.
Finally, by combining a descent estimate with
the joint stationarity bounds, we established
convergence of the directional-stationarity measures and complexity
guarantees for inexact proximal-gradient sequences.




\bibliographystyle{siamplain}
\bibliography{references}

@book{bauschke2017convex,
  author    = {Heinz H. Bauschke and Patrick L. Combettes},
  title     = {Convex Analysis and Monotone Operator Theory in Hilbert Spaces},
  series    = {CMS Books in Mathematics},
  edition   = {2nd},
  publisher = {Springer},
  address   = {Cham},
  year      = {2017},
  doi       = {10.1007/978-3-319-48311-5},
  isbn      = {978-3-319-48310-8}
}

@article{mei2026admm,
  author  = {Shenghan Mei and Chengyu Ke and Yifei Lou and Miju Ahn},
  title   = {{ADMM} and Linearized {ADMM} for Weakly Convex Minimization},
  journal = {Frontiers in Applied Mathematics and Statistics},
  volume  = {12},
  year    = {2026},
  doi     = {10.3389/fams.2026.1898659}
}

@book{rockafellar1998variational,
  author    = {R. Tyrrell Rockafellar and Roger J.-B. Wets},
  title     = {Variational Analysis},
  series    = {Grundlehren der mathematischen Wissenschaften},
  volume    = {317},
  publisher = {Springer},
  address   = {Berlin, Heidelberg},
  year      = {1998},
  doi       = {10.1007/978-3-642-02431-3},
  isbn      = {978-3-540-62772-2}
}

@misc{chen2022proximal,
  author        = {Beier Chen and Hui Zhang},
  title         = {On Proximal Gradient Mapping and Its Minimization in Norm via Potential Function-Based Acceleration},
  year          = {2022},
  eprint        = {2212.07149},
  archivePrefix = {arXiv},
  primaryClass  = {math.OC}
}

@book{beck2017firstorder,
  author    = {Amir Beck},
  title     = {First-Order Methods in Optimization},
  publisher = {Society for Industrial and Applied Mathematics},
  address   = {Philadelphia, PA},
  year      = {2017}
}

@article{davis2019stochastic,
  author  = {Damek Davis and Dmitriy Drusvyatskiy},
  title   = {Stochastic Model-Based Minimization of Weakly Convex Functions},
  journal = {SIAM Journal on Optimization},
  volume  = {29},
  number  = {1},
  pages   = {207--239},
  year    = {2019},
  doi     = {10.1137/18M1178244}
}

@article{barre2023principled,
  author  = {Mathieu Barr{\'e} and Adrien B. Taylor and Francis Bach},
  title   = {Principled Analyses and Design of First-Order Methods with Inexact Proximal Operators},
  journal = {Mathematical Programming},
  volume  = {201},
  pages   = {185--230},
  year    = {2023},
  doi     = {10.1007/s10107-022-01903-7}
}

@article{bellocruz2025relative,
  author  = {Yunier Bello-Cruz and Max L. N. Gon{\c{c}}alves and Jefferson G. Melo and Cassandra Mohr},
  title   = {A Relative Inexact Proximal Gradient Method with an Explicit Linesearch},
  journal = {Journal of Optimization Theory and Applications},
  volume  = {206},
  number  = {1},
  note    = {Article 9},
  year    = {2025},
  doi     = {10.1007/s10957-025-02684-7}
}

@article{khanh2025inexact,
  author  = {Pham Duy Khanh and Boris S. Mordukhovich and Vo Thanh Phat and Dat Ba Tran},
  title   = {Inexact Proximal Methods for Weakly Convex Functions},
  journal = {Journal of Global Optimization},
  volume  = {91},
  number  = {3},
  pages   = {611--646},
  year    = {2025},
  doi     = {10.1007/s10898-024-01460-7}
}

@article{taylor2017exact,
  author  = {Adrien B. Taylor and Julien M. Hendrickx and Fran{\c{c}}ois Glineur},
  title   = {Exact Worst-Case Performance of First-Order Methods for Composite Convex Optimization},
  journal = {SIAM Journal on Optimization},
  volume  = {27},
  number  = {3},
  pages   = {1283--1313},
  year    = {2017},
  doi     = {10.1137/16M108104X}
}

@inproceedings{rotaru2025tight,
  author    = {Teodor Rotaru and Panagiotis Patrinos and Fran{\c{c}}ois Glineur},
  title     = {Tight Analysis of Difference-of-Convex Algorithm ({DCA}) Improves Convergence Rates for Proximal Gradient Descent},
  booktitle = {Proceedings of the 28th International Conference on Artificial Intelligence and Statistics},
  series    = {Proceedings of Machine Learning Research},
  volume    = {258},
  pages     = {4114--4122},
  year      = {2025}
}

@article{rubbens2026constraint,
  author  = {Anne Rubbens and Julien M. Hendrickx},
  title   = {A Constraint-Based Approach to Function Interpolation, with Application to Performance Estimation for Weakly Convex Optimisation},
  journal = {Mathematical Programming},
  year    = {2026},
  doi     = {10.1007/s10107-026-02353-1}
}

@article{vernimmen2026worstcase,
  author  = {Pierre Vernimmen and Fran{\c{c}}ois Glineur},
  title   = {Worst-Case Convergence Analysis of Relatively Inexact Gradient Descent on Smooth Convex Functions},
  journal = {Neurocomputing},
  volume  = {660},
  pages   = {131883},
  year    = {2026},
  doi     = {10.1016/j.neucom.2025.131883}
}

@article{bach2011sparsity,
  author  = {Francis Bach and Rodolphe Jenatton
             and Julien Mairal and Guillaume Obozinski},
  title   = {Optimization with Sparsity-Inducing Penalties},
  journal = {Foundations and Trends in Machine Learning},
  volume  = {4},
  number  = {1},
  pages   = {1--106},
  year    = {2011},
  doi     = {10.1561/2200000015}
}

@article{combettes2005signal,
  author  = {Patrick L. Combettes and Val{\'e}rie R. Wajs},
  title   = {Signal Recovery by Proximal Forward-Backward Splitting},
  journal = {Multiscale Modeling and Simulation},
  volume  = {4},
  number  = {4},
  pages   = {1168--1200},
  year    = {2005},
  doi     = {10.1137/050626090}
}

@article{beck2009fista,
  author  = {Amir Beck and Marc Teboulle},
  title   = {A Fast Iterative Shrinkage-Thresholding Algorithm
             for Linear Inverse Problems},
  journal = {SIAM Journal on Imaging Sciences},
  volume  = {2},
  number  = {1},
  pages   = {183--202},
  year    = {2009},
  doi     = {10.1137/080716542}
}

@article{parikh2014proximal,
  author  = {Neal Parikh and Stephen Boyd},
  title   = {Proximal Algorithms},
  journal = {Foundations and Trends in Optimization},
  volume  = {1},
  number  = {3},
  pages   = {127--239},
  year    = {2014},
  doi     = {10.1561/2400000003}
}

@inproceedings{schmidt2011inexact,
  author    = {Mark Schmidt and Nicolas Le Roux
               and Francis Bach},
  title     = {Convergence Rates of Inexact Proximal-Gradient
               Methods for Convex Optimization},
  booktitle = {Advances in Neural Information Processing Systems},
  volume    = {24},
  pages     = {1458--1466},
  year      = {2011}
}

@article{ahn2017difference,
  author  = {Miju Ahn and Jong-Shi Pang and Jack Xin},
  title   = {Difference-of-Convex Learning: Directional Stationarity,
             Optimality, and Sparsity},
  journal = {SIAM Journal on Optimization},
  volume  = {27},
  number  = {3},
  pages   = {1637--1665},
  year    = {2017},
  doi     = {10.1137/16M1084754}
}

@article{taylor2018exactpg,
  author  = {Adrien B. Taylor and Julien M. Hendrickx
             and Fran{\c{c}}ois Glineur},
  title   = {Exact Worst-Case Convergence Rates of the
             Proximal Gradient Method for Composite
             Convex Minimization},
  journal = {Journal of Optimization Theory and Applications},
  volume  = {178},
  number  = {2},
  pages   = {455--476},
  year    = {2018},
  doi     = {10.1007/s10957-018-1298-1}
}

@article{bonettini2020forwardbackward,
  author  = {Silvia Bonettini and Marco Prato
             and Simone Rebegoldi},
  title   = {Convergence of Inexact Forward--Backward
             Algorithms Using the Forward--Backward Envelope},
  journal = {SIAM Journal on Optimization},
  volume  = {30},
  number  = {4},
  pages   = {3069--3097},
  year    = {2020},
  doi     = {10.1137/19M1254155}
}

@misc{lauga2026characterizations,
  author        = {Guillaume Lauga and Samuel Vaiter},
  title         = {Characterizations of Inexact Proximal Operators},
  year          = {2026},
  eprint        = {2602.02022},
  archivePrefix = {arXiv},
  primaryClass  = {math.OC}
}

@article{bednarczuk2026forwardbackward,
  author  = {Ewa Bednarczuk and Giovanni Bruccola
             and Gabriele Scrivanti and The Hung Tran},
  title   = {Forward--Backward Algorithms for Weakly Convex Problems},
  journal = {Applied Mathematics \& Optimization},
  volume  = {93},
  note    = {Article 74},
  year    = {2026},
  doi     = {10.1007/s00245-026-10420-4}
}

@article{hamadouche2024sharper,
  author  = {Anis Hamadouche and Yun Wu and
             Andrew M. Wallace and Jo{\~a}o F. C. Mota},
  title   = {Sharper Bounds for Proximal Gradient
             Algorithms with Errors},
  journal = {SIAM Journal on Optimization},
  volume  = {34},
  number  = {1},
  pages   = {278--305},
  year    = {2024},
  doi     = {10.1137/22M1480161}
}

@inproceedings{gu2018inexact,
  author    = {Bin Gu and De Wang and Zhouyuan Huo
               and Heng Huang},
  title     = {Inexact Proximal Gradient Methods for
               Non-Convex and Non-Smooth Optimization},
  booktitle = {Proceedings of the AAAI Conference
               on Artificial Intelligence},
  volume    = {32},
  pages     = {3093--3100},
  year      = {2018},
  doi       = {10.1609/aaai.v32i1.11802}
}
\end{document}